\RequirePackage{fix-cm}
\documentclass[smallextended]{svjour3}       % onecolumn (second format)
\smartqed  % flush right qed marks, e.g. at end of proof
\usepackage{lipsum}
\usepackage{amsfonts}
\usepackage{graphicx}
\usepackage{epstopdf}
\usepackage{algorithm}
\usepackage{algpseudocode}
\usepackage{siunitx} 
\usepackage{multirow}
\usepackage{tikz}
\usepackage{threeparttable}
\usepackage{amssymb}
\usepackage{subcaption}
\usepackage{listings}
\usepackage{xcolor}
\usepackage{amsmath}
\usepackage[hidelinks]{hyperref}
\usepackage{color,soul}
\sethlcolor{orange}
\renewcommand{\hl}[1]{#1}
\usepackage[pagewise]{lineno}

\DeclareRobustCommand{\hlcyan}[1]{{\sethlcolor{cyan}\hl{#1}}}
\DeclareRobustCommand{\hlyl}[1]{{\sethlcolor{yellow}\hl{#1}}}
\ifpdf
  \DeclareGraphicsExtensions{.eps,.pdf,.png,.jpg}
\else
  \DeclareGraphicsExtensions{.eps}
\fi

\RequirePackage{titlesec}
\titleformat{\section}  
  {\normalfont\large\bfseries}  
  {\thesection}  
  {1em}  
  {}  

\titleformat{\subsection}  
  {\normalfont\bfseries}  
  {\thesubsection}  
  {1em}  
  {}  

\titleformat{\subsubsection}  
  {\normalfont\bfseries}  
  {\thesubsubsection}  
  {1em}  
  {}  

\newcommand{\boldline}{%  
  \noalign{\smallskip}%  
  \noalign{\hrule height 1pt}%  
  \noalign{\smallskip}%  
}  

\begin{document}
%\linenumbers
\title{Matrix-Free Parallel Scalable Multilevel Deflation Preconditioning for Heterogeneous Time-Harmonic Wave Problems}


\author{Jinqiang Chen  \and Vandana Dwarka	\and  Cornelis Vuik}

\titlerunning{Parallel Matrix-Free Multilevel Deflation} % if too long for running head

\institute{
	Jinqiang Chen, corresponding author \at Delft Institute of Applied Mathematics, Delft University of Technology \\
	\email{j.chen-11@tudelft.nl}%\\
        \and
	Vandana Dwarka \at Delft Institute of Applied Mathematics, Delft University of Technology \\
	\email{v.n.s.r.dwarka@tudelft.nl}%\\
        \and
	Cornelis Vuik \at Delft Institute of Applied Mathematics, Delft University of Technology \\
	\email{c.vuik@tudelft.nl}
}

\date{Received: date / Accepted: date}

\maketitle

\begin{abstract}
We present a matrix-free parallel scalable multilevel deflation preconditioned method for heterogeneous time-harmonic wave problems. Building on the higher-order deflation preconditioning \hlcyan{proposed by Dwarka and Vuik (SIAM J. Sci. Comput. 42(2):A901–A928, 2020; J. Comput. Phys. 469:111327, 2022)} for highly indefinite time-harmonic waves, we adapt these techniques for parallel implementation in the context of solving large-scale heterogeneous problems with minimal pollution error. \hlcyan{Our proposed method integrates the Complex Shifted Laplacian preconditioner with deflation approaches.} We employ higher-order deflation vectors and re-discretization schemes derived from the Galerkin coarsening approach for a matrix-free parallel implementation. We suggest a robust and efficient configuration of the matrix-free multilevel deflation method, which yields a close to wavenumber-independent convergence and good time efficiency. Numerical experiments demonstrate the effectiveness of our approach for increasingly complex model problems. The matrix-free implementation of the preconditioned Krylov subspace methods reduces memory consumption, and the parallel framework exhibits satisfactory parallel performance and weak parallel scalability. This work represents a significant step towards developing efficient, scalable, and parallel multilevel deflation preconditioning methods for large-scale real-world applications in wave propagation.
\end{abstract}

\section{Introduction}
The Helmholtz equation is a crucial mathematical model that describes the behavior of time-harmonic waves in various scientific fields, such as seismology, sonar technology, and medical imaging. \hlcyan{While the classical formulation with standard Laplacian operator remains essential for many practical applications, researchers have also explored extended formulations to address specific physical phenomena. Recent studies have investigated Helmholtz equations in various forms, including elastic Helmholtz equations {\cite{yovelLFAtunedMatrixfreeMultigrid2024}}, the Stochastic Helmholtz equation {\cite{pulchHelmholtzEquationUncertainties2024}}, and Helmholtz equations with fractional Laplacian operators \mbox{\cite{liPreconditioningTechniqueBased2023,adrianiAsymptoticSpectralProperties2024}}. In this work, we focus on the classical acoustics Helmholtz equation.} Solving this equation numerically involves dealing with a sparse, symmetric, complex-valued, non-Hermitian, and indefinite linear system. For large-scale problems, iterative methods and parallel computing are commonly used. However, the indefiniteness of the system introduces significant difficulties, particularly when high frequencies are involved, as it limits the convergence of iterative solvers. Furthermore, to control pollution errors, it is essential to refine the grid so that $k^3h^2 < 1$ \cite{babuska1997pollution}, where $k$ is the wavenumber and $h$ is the mesh size. The challenge of efficiently solving the Helmholtz equation while preserving high accuracy and minimizing pollution errors continues to be an active area of research. \hlcyan{The development of a parallel scalable iterative method with convergence properties independent of the wavenumber could have far-reaching implications for various disciplines, including electromagnetics, seismology, and acoustics \mbox{\cite{sourbierThreedimensionalParallelFrequencydomain2011,tournierNumericalModelingHighspeed2017}}}.

A variety of preconditioners have been proposed for the Helmholtz problem, and among them, the Complex Shifted Laplace Preconditioner (CSLP) \cite{erlangga2004class,erlangga2006novel} is one of the most popular options in the industry. The CSLP exhibits good properties for medium wavenumbers. However, the eigenvalues shift towards the origin as the wavenumber increases. As a result, the deflation method was introduced to accelerate the convergence of the CSLP-preconditioned Krylov subspace method \cite{erlangga2008multilevel,sheikh2016accelerating}. However, the number of iterations in both variations still gradually increases with the wavenumber $k$. In a recent development, Dwarka and Vuik \cite{dwarka2020scalable} introduced higher-order approximation schemes to construct deflation vectors. This two-level deflation method exhibits convergence that is nearly independent of the wavenumber. The authors further extend the two-level deflation method to a multilevel deflation method \cite{dwarka2020scalablemultilevel}. Using higher-order deflation vectors, the authors demonstrate that the near-zero eigenvalues of the coarse-grid operators remain aligned with those of the fine-grid operator up to the level where the coarse-grid linear systems become negative indefinite. This alignment prevents the spectrum of the preconditioned system from approaching the origin. By combining this approach with the CSLP preconditioner, the authors achieved an iterative method with close to wavenumber-independent convergence for highly indefinite linear systems. Incorporating the deflation preconditioner has resulted in improved convergence; however, it has an impact on efficiency in relation to memory and computational cost. \hlyl{A promising branch is the use of Domain Decomposition Methods (DDM) as preconditioning techniques. These methods typically require two essential components: carefully designed transmission conditions and problem-adapted coarse spaces. Notable developments include the DtN and GenEO spectral coarse spaces {\cite{bootlandComparisonCoarseSpaces2021}}, which utilize selected modes from local eigenvalue problems specifically tailored to the Helmholtz equation. For comprehensive surveys on DDM preconditioners for the Helmholtz problems, we refer the reader to {\cite{gander2019class}} and the references therein. An alternative direction is the Multiscale Generalized Finite Element Method (MS-GFEM) {\mbox{\cite{babuskaOptimalLocalApproximation2011,maNovelDesignAnalysis2022}}}. Recent work by Ma et al. {\cite{chupengWavenumberExplicitConvergence2023}} has successfully applied MS-GFEM with novel local approximation spaces to high-frequency heterogeneous Helmholtz problems. While the deflation method achieves near wavenumber-independent convergence by aligning the near-zero eigenvalues between coarse and fine-grid operators, the discrete MS-GFEM approach solves the problem in one shot without iterating based on solving some carefully-designed local problems and a global coarse problem, suggesting potential benefits in combining both approaches for future research.}

Efforts are also underway to develop parallel scalable Helmholtz solvers. With well-designed parallelization strategies, domain decomposition methods have shown promise in reducing the number of iterations and improving computational efficiency \cite{taus2020sweeps}. Another approach is the parallelization of existing advanced algorithms. Parallel implementations of Bi-CGSTAB preconditioned by multigrid-based CSLP have been presented for 2D and 3D forward modeling by Kononov and Riyanti \cite{Kononov2007Numerical,riyanti2007parallel}, respectively. Gordon and Gordon \cite{Gordon2013Robust} proposed the block-wise parallel extension of their so-called CARP-CG algorithm (Conjugate Gradient Acceleration of CARP). The block-parallel CARP-CG algorithm shows improved scalability as the wavenumber increases. Calandra et al. \cite{calandra2013improved,calandra2017geometric} proposed a geometric two-grid preconditioner for 3D Helmholtz problems, which exhibits strong scaling in massively parallel setups. 

Although conventional multigrid methods using standard smoothing and coarse grid corrections fail for the Helmholtz equation, their high efficiency in solving positive definite problems has motivated research into developing robust multilevel approaches for this equation \cite{Elman2001,kim2002multigrid,lu2019robust,dwarka2020scalablemultilevel}. \hlcyan{While previous works have established the theoretical foundations for multilevel deflation methods \mbox{\cite{tang2009comparison,sheikh2013convergence,sheikh2016accelerating,dwarka2020scalablemultilevel}}, our focus is on a parallel scalable implementation of multilevel deflation for practical large-scale applications. We aim to perform comprehensive numerical experiments to validate the theoretical predictions and demonstrate the method's effectiveness in large-scale scenarios, where parallel implementation challenges often exceed idealized theoretical assumptions.} 

Based on the parallel framework of the CSLP-preconditioned Krylov subspace methods \cite{jchen2D2022,etna_vol59_pp270-294}, a matrix-free parallel two-level deflation preconditioning \cite{chen2023matrixfree2ldef} has been implemented recently. To the author's knowledge, there is no existing literature on parallel multilevel deflation so far. This paper addresses this gap by proposing a matrix-free, parallel scalable multilevel deflation preconditioning method. In this work, we explore methods to extend the wavenumber-independent convergence from the two-level to a multilevel setting. \hlcyan{This work presents significant innovations in solving large-scale Helmholtz problems. We develop novel re-discretization schemes for multilevel hierarchies, ensuring effective approximation of Galerkin coarsening operators across all levels while maintaining the matrix-free parallel framework. Through comprehensive numerical experiments, we establish optimal parameters for robust convergence across different problem scales. Furthermore, we introduce a controllable tolerance for coarse-level iterations, a previously unexplored but essential component for achieving wavenumber-independent convergence in practical multilevel deflation methods. These innovations culminate in a highly efficient parallel framework that demonstrates both wavenumber-independent convergence and excellent scaling properties in massively parallel environments, as validated by extensive numerical experiments.}

% The outline is not required, but we show an example here.
The paper is organized as follows. In Section 2, we begin by describing our model problems. We present the matrix-free parallel variant of the multilevel deflated Krylov method in Section 3. Section 4 presents an optimally tuned configuration of the matrix-free parallel multilevel deflation method. Finally, we present numerical results to evaluate parallel performance in Section 5. Section 6 contains our conclusions.

\section{Problem description}
\label{sec:MPs}

We mainly consider the following two-dimensional mathematical model. Suppose that the domain $\Omega$ is rectangular with a boundary $\Gamma = \partial\Omega$. The Helmholtz equation reads
\begin{equation}\label{eq:HelmholtzEq}
-\Delta \mathbf{u} - k(x,y)^2 \mathbf{u} = \mathbf{b},\  \text{on} \ \Omega, \nonumber
\end{equation}
supplied with either Dirichlet or Sommerfeld radiation boundary conditions. Suppose the frequency is $f$, the speed of propagation is $c(x,y)$, they are related by
\begin{equation}\label{eq:wavenum}
k(x,y) =\frac{2\pi f}{c(x,y)}. \nonumber
\end{equation}

% one of the following boundary conditions:
% \begin{equation}\label{eq:DirichletBC}
% \text{Dirichlet: } \mathbf{u}=\mathbf{g},\  \text{on} \ \partial\Omega, \nonumber
% \end{equation}
% \begin{equation}\label{eq:SommerfeldBC}
% \text{first-order Sommerfeld: } \frac{\partial \mathbf{u}}{\partial \vec{n}}-\text{i} k(x,y) \mathbf{u} = \mathbf{0},\  \text{on} \ \partial\Omega, \nonumber 
% \end{equation}
% where $\text{i}$ is the imaginary unit. $\vec{n}$ and $\mathbf{g}$ represent the outward normal and the given data of the boundary, respectively. $\mathbf{b}$ is the source function. $k(x,y)$ is the wavenumber on $\Omega$.

\subsection{Discretization} \label{sec:Discretization}
The computational domain is discretized using structural vertex-centered grids with uniform mesh width $h$. The discrete approximation of $u(x,y)$ is denoted as $u(i,j)$ or $u_{i,j}$, where grid points $(x_i, y_j)$ are given by $x_i = x_1 + (i-1)h$ and $y_j = y_1 + (j-1)h$.

A second-order finite difference scheme for a 2D Laplace operator has a well-known $3\times3$ stencil.
% \begin{equation}\label{eq:stencilLaplace}
% \left[ {{-\Delta_h}} \right] = \frac{1}{{{h^2}}}\left[ {\begin{array}{*{20}{c}}
%     0&{ - 1}&0\\
%     { - 1}&{4}&{ - 1}\\
%     0&{ - 1}&0
%     \end{array}} \right]. \nonumber
% \end{equation}
Similarly, we denote a computation stencil for the wavenumber term in the Helmholtz equation as
\begin{equation}\label{eq:stencilwavenumber}
\left[ \mathcal{I}(k_{i,j}^2)_h \right] = \left[ {\begin{array}{*{20}{c}}
    0&0&0\\
    0&k_{i,j}^2&0\\
    0&0&0
    \end{array}} \right], \nonumber
\end{equation}
where $\mathcal{I}(k_{i,j}^2)$ represents a diagonal matrix with $k_{i,j}^2$ as its diagonal elements. The stencil of the Helmholtz operator $A_h$ can be obtained by
% subtracting the diagonal matrix $ \mathcal{I}(k_{i,j}^2)_h$  to the Laplacian operator $-\Delta_h$, \textit{i.e.}
\begin{equation}\label{eq:helm operator h}
\left[A_h\right] =\left[-\Delta_h\right]- \left[\mathcal{I}(k_{i,j}^2)_h\right].
\end{equation}
% Therefore, the stencil of the discrete Helmholtz operator is
% \begin{equation}\label{eq:Helmstencil}
% \left[ {{A_h}} \right] = \frac{1}{{{h^2}}}\left[ {\begin{array}{*{20}{c}}
%     0&{ - 1}&0\\
%     { - 1}&{4 - {k_{i,j}^2}{h^2}}&{ - 1}\\
%     0&{ - 1}&0
%     \end{array}} \right].
% \end{equation}

For boundary conditions, a ghost point located outside the boundary points can be introduced.
For instance, suppose $u_{0,j}$ is a ghost point on the left of $u_{1,j}$, 
% the normal derivative can be approximated  by
% \begin{equation}
%     \label{eq:GhostPoint}
%     \frac{\partial u}{\partial \vec{n}} -\text{i} k(x,y) u =\frac{u_{0,j}-u_{2,j}}{2h}-\text{i} k_{1,j} u_{1,j}= 0. \nonumber
% \end{equation}
for Sommerfeld radiation boundary condition, we have
\begin{equation}\label{eq:ghost_sommerfeld}
    u_{0,j} = u_{2,j} + 2h\text{i} k_{1,j} u_{1,j}.
\end{equation}
For the Dirichlet boundary condition, we have
\begin{equation}\label{eq:ghost_Dirichlt}
    u_{0,j} = 2u_{1,j} - u_{2,j}.
\end{equation}

Discretization of the partial equation on the finite-difference grids results in a system of linear equations $A_h \mathbf{u}_h=\mathbf{b}_h$. With first-order Sommerfeld radiation boundary conditions, the resulting matrix is sparse, symmetric, complex-valued, indefinite, and non-Hermitian.

Note that $kh$ is an important parameter that indicates how many grid points per wavelength are needed. The mesh width  $h$ can be determined by the guidepost of including at least $N_{pw}$(e.g. 10 or 30) grid points per wavelength. They have the following relationships 
\begin{equation}
kh=\frac{2 \pi h}{\lambda} = \frac{2 \pi }{N_{pw}}. \nonumber
\end{equation}
For example, if at least 10 grid points per wavelength are required, we can maintain $kh = 0.625$.

\subsection{Model Problem - constant wavenumber}
We first consider a 2D problem with constant wavenumber in a rectangular homogeneous domain $\Omega=\left[0,1 \right]$. 
% A point source is given by
% \begin{equation}\label{eq: dirac rhs}
% b\left( {x,y} \right) = \delta \left( {x - x_0,y - y_0} \right),\ \Omega=\left[0,1 \right],  \nonumber
% \end{equation}
% where $\delta \left( {x,y} \right)$ is a Dirac delta function in the following form
% \begin{equation}\label{eq: dirac func}
% \delta \left( {x,y} \right) = \left\{ {\begin{array}{*{20}{c}}
%     { + \infty \;\;\;x = 0,\;y = 0}\\
%     {0\;\;\;\;\;\;\;x \ne 0,\;y \ne 0}
%     \end{array}} \right., \nonumber
% \end{equation}
% which satisfies
% \begin{equation}\label{eq: dirac func restricton}
% \int {\int {{\delta ^2}\left( {x,y} \right)dxdy = 1} }. \nonumber
% \end{equation}
A point source defined by a Dirac delta function is imposed at the center $(x_0,y_0)=(0.5,0.5)$. The wave propagates outward from the center of the domain. The Dirichlet boundary conditions (denoted as MP-1a) or the first-order Sommerfeld radiation boundary conditions (denoted as MP-1b) are imposed, respectively.  

\subsection{Model Problem - Wedge problem}
\begin{figure}[htbp]
    \centering
    \includegraphics[width=0.38\textwidth, trim=10 10 10 45,clip]{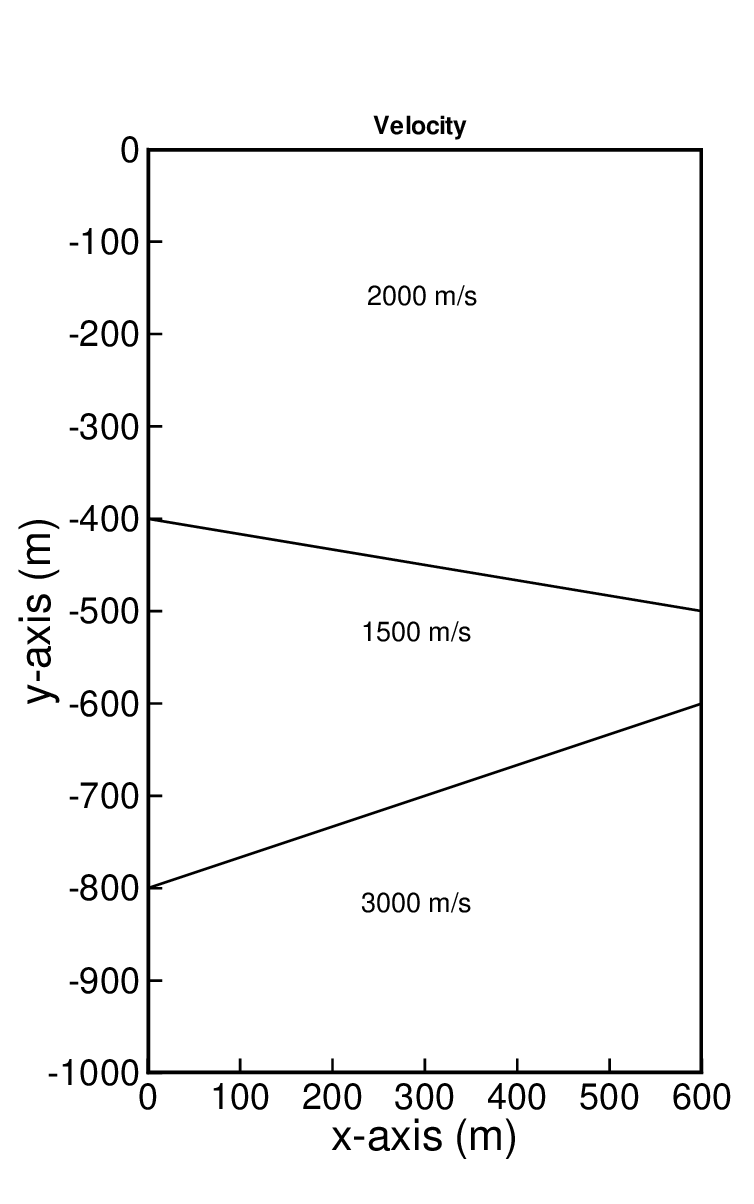}
    \caption{The velocity distribution of the Wedge problem}
    \label{fig:wedge domain}
\end{figure}

Most physical problems of geophysical seismic imaging describe a heterogeneous medium. The so-called Wedge problem \cite{plessix2003separation} is a typical problem with simple heterogeneity. It mimics three layers with different velocities, hence different wavenumbers. As shown in Figure \ref{fig:wedge domain}, the rectangular domain $\Omega=\left[0,600 \right] \times \left[-1000,0 \right]$  is split into three layers. Suppose the wave velocity $c$ is constant within each layer but different from each other. A point source is located at $\left(x, y \right) = \left( 300, 0\right) $. The wave velocity $c(x,y)$ is shown in Figure \ref{fig:wedge domain}. The first-order Sommerfeld radiation boundary conditions are imposed on all boundaries.

\subsection{Model Problem - Marmousi Problem} 
\label{sec:Marmousi}
\begin{figure}[htbp]
    \centering
    \includegraphics[width=0.85\textwidth, trim=4 4 4 4,clip]{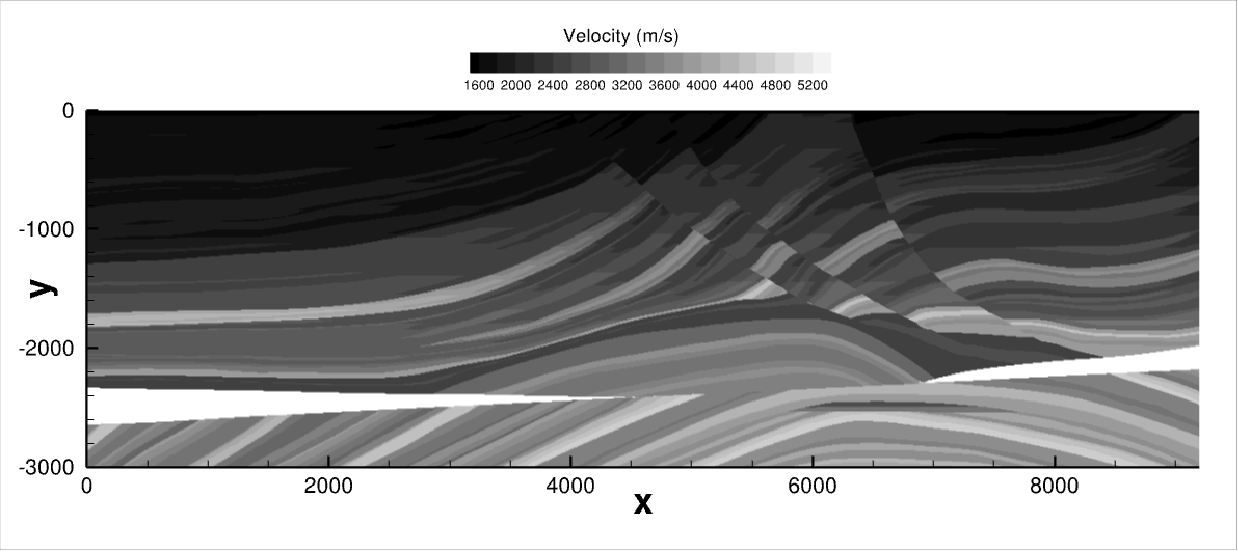}
    \caption{The velocity distribution  of the Marmousi problem}
    \label{fig:marmousi}
\end{figure}

For industrial applications, the third model problem is the so-called Marmousi problem \cite{Versteeg_1991_ME}, a well-known benchmark problem. It contains 158 horizontal layers in the depth direction, making it highly heterogeneous. The wave velocity $c(x,y)$ over the domain is shown in Figure \ref{fig:marmousi}. The first-order Sommerfeld radiation boundary conditions are imposed on all boundaries.

\section{Deflated Krylov methods}
\label{sec:NMs}
We review the two-level deflation preconditioning and a preliminary multilevel setting, then detail their adaptation to a matrix-free parallel framework.

\subsection{Two-level deflation}
Suppose a general nonsingular linear system $Au=b$, where $A\in\mathbb{R}^{n \times n}$, and a projection subspace matrix, $Z \in \mathbb{R}^{n\times m}$, with $m<n$ and full rank are given. Assume that $E=Z^T A Z$ is invertible, the projection matrix $P$ can be defined as
\begin{equation}
    \label{eq: def DEF}
    P = I - AQ, \quad Q=Z E^{-1} Z^T, \quad E = Z^T A Z.
\end{equation}
The observation that multigrid inter-grid operators emphasize small frequencies and preserve them on coarser levels leads to the possibility of using geometrically constructed multigrid vectors as deflation vectors. In such scenarios, we refer to $E$ as the coarse grid Helmholtz operator, which exhibits similar properties to that of $A$. 

The CSLP-preconditioner can be included to obtain even better convergence. The CSLP preconditioner $M_{(\beta_1,\beta_2)}$ is defined by
\begin{equation}\label{eq: CSLP def}
    M_{(\beta_1,\beta_2)}=-\Delta-(\beta_1+\text{i}\beta_2) \mathcal{I}(k_{i,j}^2),
\end{equation}
where $\text{i}=\sqrt{-1}$, $(\beta_1,\beta_2)\in [0,1]$. \hl{The solution of complex-shifted Laplacian systems, which arise from approximating the inverse of CSLP, has been extensively studied in the literature. Various approaches have been proposed, including multigrid methods {\mbox{\cite{erlangga2006novel,cocquetHowLargeShift2017}}} and Krylov subspace methods {\mbox{\cite{grahamDomainDecompositionPreconditioning2017,dwarka2020scalablemultilevel}}}.  A recent alternative approach by {\cite{linAbsolutevalueBasedPreconditioner2024}} introduces absolute-value based preconditioners for the equivalent block real linear system formulation of the complex-shifted Laplacian system.} The multigrid-based CSLP mentioned in this paper will fully adopt the matrix-free parallel framework and settings proposed in \cite{jchen2D2022}.

To allow approximate solvers for $E^{-1}$, one can prevent close-to-zero eigenvalues from obstructing the convergence of the Krylov solver by adding a term $Q$ to deflate toward the largest eigenvalues of the preconditioned system \cite{tang2009comparison}. As a result, the Adapted Deflation Variant 1 (A-DEF1) preconditioner $P$ reads as
\begin{equation}
    P := M^{-1}_{(\beta_1,\beta_2)}P+Q. \nonumber
\end{equation}
When higher-order Bezier curves are used to construct high-order deflation vectors $Z$, this results in the so-called Adapted Preconditined Deflation (APD) \cite{dwarka2020scalable}. The preconditioned linear system to be solved becomes
\begin{equation}
\label{eq: ADP linear sys}
    PAu = Pb.
\end{equation}
Note that $PA$ is nonsingular, so Eq. (\ref{eq: ADP linear sys}) has a unique solution. In this work, we will stick to the use of the higher-order deflation vectors.

\subsection{Multilevel deflation}
When using the two-level method in practical large-scale applications, solving the coarse-grid system remains expensive, whether solved exactly \cite{dwarka2020scalablemultilevel}, or approximately by CSLP-preconditioned Krylov methods \cite{chen2023matrixfree2ldef}. In accordance with the multigrid method, as the coarse grid system has similar properties as the original Helmholtz operator, one can obtain a multilevel framework by applying the two-level cycle recursively, as shown in Algorithm \ref{algthm: recursively TL DEF}. The flexible subspace Krylov method such as FGMRES preconditioned by two-level deflation is applied recursively on subsequent coarse-grid systems $E$. Compared to the multilevel deflation proposed in \cite{dwarka2020scalablemultilevel}, a few remarks are noted here.

\begin{algorithm}[htbp]
    \caption{Recursive two-level deflated FGMRES: \texttt{TLADP-FGMRES(A, b)}} \label{algthm: recursively TL DEF}
    \begin{algorithmic}[1]
        \State Determine the current level $l$ and dimension $m$ of the Krylov subspace
        \State{Initialize $u_0$, compute {$r_0 =b-Au_0$}, $\beta=\left| \left| r_0\right| \right|$, $v_1 = r_0/\beta $;}
        \State{Define $\bar{H}_m \in \mathbb{C}^{(m+1) \times m}$ and initialize to zero}
        \For{$j=1,2,...,m$ or until convergence}
        \State{$\hat{v}_j=Z^T v_j$ } \Comment{Restriction}
        \If{\hl{$l+1 == ml$}} \Comment{\hl{Predefined coarsest level $ml$}}  
        \State $\tilde{v} \approx E^{-1}\hat{v}$ \Comment{Approximated by CSLP-FGMRES}
        \Else
        \State{$l \gets l+1$}
        \State{$\tilde{v} \gets $\texttt{TLADP-FGMRES(E, $\hat{v}$)} } \Comment{Apply two-level deflation recursively}
        \EndIf
        \State{$t=Z\tilde{v}$} \Comment{Interpolation}
        \State{$s=At$}
        \State{$\tilde{r}=v_j-s$}
        \State{$r \approx M^{-1}\tilde{r}$ \Comment{CSLP, by multigrid method or Krylov iterations}}
        \State{$x_j=r+t$}  
        \State{$w =A x_j$ }
        %\State{\color{red}$w=M^{-1} \tilde{v}_j$ \ \% Preconditioned}
        \For{$i:=1,2,...,j$}
        \State{$h_{i,j}=\left( w, v_i\right) $}
        \State{$w \gets w-h_{i,j}v_i$}
        \EndFor
        \State{$h_{j+1,j}:=\left| \left| w \right| \right|_2$, $v_{j+1}=w/h_{j+1,j}$}
        \EndFor
        \State{\hl{$X_m=\left[ x_1,...,x_m \right]$, $\bar{H}_m=\left\lbrace h_{i,j}\right\rbrace_{1\le i \le j+1, 1\le j \le m}$}}
        \State{$u_m=u_0+X_my_m$ where $y_m = \text{arg min}_y\left| \left| \beta e_1 - \bar{H}_m y\right| \right|$}
        \State{\textbf{Return} $u_m$}
    \end{algorithmic}
\end{algorithm}

First, Algorithm \ref{algthm: recursively TL DEF} does not include the process of determining the corresponding coarser-grid system $E$ and CSLP preconditioner $M$ on the current level $l$. This will be elaborated on in the next subsection.

Second, for efficient parallelization, we employ a GMRES method to solve the coarsest grid problem approximately rather than a direct solver. This method is preconditioned by CSLP, which is defined according to the coarsest grid operator. We will numerically investigate the necessary accuracy (or the number of iterations) for solving the coarsest grid problem in the next section.

Third, \hl{the multilevel deflation method requires approximating the inverse of CSLP on each level. In {\cite{dwarka2020scalablemultilevel}}, this is accomplished using several Krylov subspace iterations (e.g., Bi-CGSTAB) on all levels.} The authors set the maximum number of iterations to $C_{it} (N^{l})^{\frac{1}{4}}$, where $C_{it}$ is a constant and $N^{l}$ denotes the problem size on level $l$. This strategy allows the benefits of using a small shift, resulting in a preconditioner similar to the original Helmholtz operator that retains the ability to shift indefiniteness at certain levels. However, the maximum number of iterations is positively correlated with the grid size on each level, indicating that larger grid sizes require more iterations. \hl{Considering the large-scale applications, utilizing Krylov subspace iterations on the first level (finest grid) or the second level may become computationally intensive. Therefore, we propose employing a multigrid cycle to approximate CSLP on the first or second level.} Several Krylov subspace iterations can then be applied on the coarser levels. However, in the case of multigrid-based CSLP, ensuring a sufficiently large complex shift is essential. In addition to setting the maximum number of iterations, a relative tolerance as stopping criteria for the iterations is also established in this paper. This allows iterations to cease once the maximum number of iterations or the tolerance is reached. 

Fourth, as shown in Algorithm \ref{algthm: recursively TL DEF}, the number of deflated FGMRES iterations is specified by $m$. The cycle type of the multilevel deflation technique is determined by the number of iterations of the deflated FGMRES on each coarse level, except for the finest level. If only one iteration is allowed on the coarser levels, a V-cycle is obtained, which is similar to the V-cycle structure of multigrid when $\gamma=1$. Correspondingly, two iterations on the coarser levels will result in a W-cycle. According to the multigrid method, the W-cycle may offer faster convergence than V-cycle but at the expense of computational efficiency. For increasingly complex model problems, striking a balance between optimal convergence and computational efficiency in the selection of $m$, hence determining the necessary accuracy (or the number of iterations) for the coarser levels, will be a focal point of this study.

Fifth, in Algorithm \ref{algthm: recursively TL DEF}, all involved matrix-vector multiplications (lines 5, 7, 10, 12, 13, 15, 17) are expressed to denote the outcome of these operations. In our implementation, we compute and return the result of matrix-vector multiplication based on input variables through a linear combination. No explicit construction of any matrices takes place in our approach.

\subsection{Matrix-free parallel implementation}
 Matrix-free implementations offer a compelling alternative to standard sparse matrix data formats in large-scale computational scenarios. Besides the reduced memory consumption, matrix-free methods exhibit performance advantages and can potentially outperform matrix-vector multiplications with stored matrices \cite{drzisga2020stencil}. \hl{Through roofline model analysis (detailed in Appendix {\ref{Appendix B}}), we demonstrate that our approach achieves 2.35 times higher arithmetic intensity compared to traditional CSR matrix-based implementations. This theoretical advantage translates to substantial performance gains, particularly for large-scale problems.} These improvements  enable the solution of larger-scale Helmholtz problems previously constrained by memory limitations, while also enhancing the applicability of modern data-driven methods \cite{drzisga2023semi}. This section details the matrix-free implementation of operators in the multilevel deflation method.
 
 Matrix-free matrix-vector multiplication is implemented using stencil notation. The computational stencils for both the finest-level and second-level operators (Helmholtz and CSLP preconditioner) and grid-transfer operators (higher-order interpolation and restriction) are detailed in \cite{chen2023matrixfree2ldef}. To enable true multilevel deflation, we extend the matrix-free implementation to coarser levels. We denote the Helmholtz operators as $A_{2^{l-1}h}$ and the CSLP operators as $M_{2^{l-1}h}$, where $l$ is the level number ($l=1$ represents the finest grid). Starting from the second-level grid, we want to find the computational stencils for $A_{4h}$ so that it is a good approximation to the Galerkin coarsening operator $Z^TA_{2h}Z$. Following \cite{chen2023matrixfree2ldef}, we decompose the Helmholtz operator into Laplace and wavenumber operators (assuming a constant wavenumber). By applying Galerkin coarsening operations to their stencils, we obtain the following stencils of the Laplace and wavenumber operators for interior points on the third-level coarse grid:
\begin{equation}
\label{eq: Glk Laplace 4h}
    \begin{aligned}
        & [-\Delta_{4h}] = \frac{1}{4096}\cdot \frac{1}{1024}\cdot \frac{1}{h^2}\cdot \\
        & \left[ {\begin{array}{*{20}{r}}
                 -3&   -534&  -5773& -11956&  -5773&   -534&    -3\\
               -534& -32844&-207370&-354088&-207370& -32844&  -534\\
              -5773&-207370&-294371& 384244&-294371&-207370& -5773\\
             -11956&-354088& 384244&2945488& 384244&-354088&-11956\\
              -5773&-207370&-294371& 384244&-294371&-207370& -5773\\
               -534& -32844&-207370&-354088&-207370& -32844&  -534\\
                 -3&   -534&  -5773& -11956&  -5773&   -534&    -3
            \end{array}} \right]_{4h},
    \end{aligned} \nonumber
\end{equation}

\begin{equation}
\label{eq: ReD-Glk K 4h}
    \begin{aligned}
        & [ \mathcal{I}(k^2)_{4h}] =\frac{1}{4096}\cdot \frac{1}{4096}\cdot k^2 \cdot \\
        & \left[ {\begin{array}{*{20}{r}}
                1&    322&    3823&    8092&    3823&    322&   1\\
              322& 103684& 1231006& 2605624& 1231006& 103684& 322\\
             3823&1231006&14615329&30935716&14615329&1231006&3823\\
             8092&2605624&30935716&65480464&30935716&2605624&8092\\
             3823&1231006&14615329&30935716&14615329&1231006&3823\\
              322& 103684& 1231006& 2605624& 1231006& 103684& 322\\
                1&    322&    3823&    8092&    3823&    322&   1
            \end{array}} \right]_{4h}.
    \end{aligned} \nonumber
\end{equation}
Using these stencils, the Helmholtz operator and CSLP operator on the third level can be obtained according to their definitions Eqs. (\ref{eq:helm operator h}) and (\ref{eq: CSLP def}), respectively. 

Continuing this process iteratively, we can obtain stencils for the Helmholtz operator on coarser levels. It should be noted that starting from the third level, the size of the computation stencils will remain at $7\times7$. Specific stencils for the fourth to sixth levels can be found in Appendix \ref{Appendix A}.

\paragraph{Boundary} 
Introducing an accurate boundary scheme for the aforementioned $7 \times 7$ computational stencils remains an open problem. In this paper, we present a simple yet effective approach, involving the introduction of a ghost point outside the physical boundaries, as depicted in Eqs. (\ref{eq:ghost_sommerfeld}) and (\ref{eq:ghost_Dirichlt}). We apply standard second-order finite-difference discretization to points on the physical boundary. For points near the boundary, we set additional grid points beyond the ghost point to zero. It is important to note that the wavenumbers of the ghost points are also required. For Dirichlet boundary conditions, we determine the wavenumbers of the ghost points similar to Eq. (\ref{eq:ghost_Dirichlt}). In other cases, the wavenumbers of the ghost points are uniformly set to zero. This zero-padding approach is motivated by the observation that the coefficients outside the $3\times3$ kernel become small, and thus, the influence of these points on the overall solution is expected to be minimal. By setting these points to zero, we aim to simplify the computation while maintaining the accuracy of the solution. 

To develop a parallel scalable iterative solver, the matrix-free multilevel deflated Krylov subspace methods are implemented within the parallel framework presented by \cite{jchen2D2022,chen2023matrixfree2ldef}. 

\section{Configuration}
Before presenting the performance analysis of the matrix-free multilevel deflation method, we systematically tune the essential components of the algorithm to achieve an optimal balance between computational efficiency and numerical robustness. This section establishes the precise configuration that ensures wavenumber-independent convergence while minimizing computational overhead for complex numerical applications. The outer FGMRES iterations start with a zero initial guess and terminate when the relative residual in Euclidean norm reaches $10^{-6}$.
Note that all presented results in this section are obtained from sequential computations. \hl{In our notation, L$n$ represents the $n$-th level in the multigrid hierarchy, where L1 corresponds to the finest level.}

\subsection{Tolerance for solving the coarsest problem} \label{sec:coarsest_tol}
In this subsection, we explore the tolerance considerations for solving the coarsest problem. For better comparison and fewer other influencing factors, we perform a V-cycle three-level deflation approach to solve the constant wavenumber model problem with Sommerfeld radiation boundary conditions. The finest level represents the first level, and the third level corresponds to the coarsest problem, which will be addressed using GMRES preconditioned with CSLP. To ensure an accurate inverse of CSLP on each level, we employ Bi-CGSTAB iterations, reducing the relative residual to $10^{-8}$, assuming that the optimal tolerance is unknown. Since the Krylov subspace method instead of the multigrid method is employed to solve the CSLP here, a small complex shift can be used. Here we will use $\beta_2=0.1$. The model problem with wavenumber $k=200$ is solved with two kinds of resolution, that is $kh=0.3125$ and $0.625$, respectively. As analyzed in \cite{dwarka2020scalablemultilevel}, the third level will remain indefinite for $kh=0.3125$, while it becomes negative definite for $kh=0.625$.

Table \ref{tab:on_tol_coarsest_problem} illustrates the impact of varying tolerances for solving the coarsest problem on the convergence behavior. Specifically, it presents the number of outer iterations required to reduce the relative residual to $10^{-6}$ and the corresponding number of iterations needed to solve the coarsest problem once.

\begin{table}[htbp]
    \centering
    \caption{The impact of varying tolerances for solving the coarsest problem on the number of outer iterations. In parentheses is the number of iterations required to solve the corresponding coarsest problem once.}
    \label{tab:on_tol_coarsest_problem}
    \scalebox{0.8}{
    \begin{tabular}{lllllll}
        \hline\noalign{\smallskip}
         $kh$ & $10^0$ & $10^{-1}$ & $10^{-2}$ & $10^{-4}$ & $10^{-6}$ & $10^{-8}$ \\
        \boldline
        $0.3125$ & 12 (1) & 5 (8) & 5 (17) & 5 (31) & 5 (44) & 5 (58) \\
        $0.625$ & 16 (1) & 31 (32) & 33 (67) & 33 (133) & 33 (203) & 33 (256) \\
        \noalign{\smallskip}\hline
    \end{tabular}}
\end{table}

From Table \ref{tab:on_tol_coarsest_problem}, we observe varying accuracy requirements for solving the coarsest grid problem corresponding to different values of $kh$. In the case of $kh=0.3125$, a relative tolerance of $10^{-1}$ is necessary for maintaining the convergence of the outer iterations. Conversely, for $kh=0.625$, a single iteration of the coarsest grid solver is sufficient. A strict tolerance in solving the coarsest grid even leads to more outer iterations.

We attribute this phenomenon to the nature of the coarsest grid system, whether it is indefinite or negative definite. According to \cite{dwarka2020scalablemultilevel}, the third level remains indefinite for $kh=0.3125$, while it becomes negative definite for $0.625$. If the coarsest grid system is indefinite, the relative tolerance of the iterative solver should be ensured at $10^{-1}$ or smaller. On the contrary, if the coarsest grid system is negative definite, one iteration is adequate. However, it leads to more outer iterations compared to the former case.

To further validate this conclusion, \hl{we note the following numerical observations (not presented in tables or figures): using the model problem with Dirichlet boundary conditions at $kh=0.625$, a tolerance of $10^0$ results in outer iterations of $27$, while $10^{-1}$ leads to $32$. This brief numerical observation supports our finding that only one iteration suffices when the system becomes negative definite on the coarsest grid.} It should be noted that, for the sake of uniformity, if a tolerance of $10^{0}$ appears in this paper, it means that only one iteration is performed.

\subsection{Tolerance for solving CSLP} \label{sec:tol_for_cslp}
In this section, we explore the accuracy requirements for the approximate inverse of CSLP on each grid level. Consistent with the solver settings from the previous subsection, we perform one iteration on the coarsest level for $kh=0.625$ and set the tolerance for the coarsest grid solver to $10^{-1}$ for $kh=0.3125$. We vary the tolerance for the Bi-CGSTAB solver used in the approximate CSLP solution. The results are presented in Table \ref{tab:tolerance-study-625} and Table \ref{tab:tolerance-study-3125}, where ``$Ln$ Bi-CGSTAB'' denotes the number of Bi-CGSTAB iterations needed to achieve the corresponding tolerance on the $n$-th level. The phrases ``Outer FGMRES'' and ``Coarsest FGMRES'' denote the number of FGMRES iterations required for the outer solvers and the coarsest problem solver, respectively.

\begin{table}[htbp]
    \centering
    \caption{Tolerance study for CSLP approximation for $kh=0.625$.}
    \label{tab:tolerance-study-625}
    \scalebox{0.8}{\begin{tabular}{cccccc}
        \hline\noalign{\smallskip}
                     & $10^{-1}$ & $10^{-2}$ & $10^{-4}$ & $10^{-6}$ & $10^{-8}$ \\ \boldline
        Outer FGMRES & 18 & 16 & 16 & 16 & 16 \\
        Coarsest FGMRES & 1 & 1 & 1 & 1 & 1 \\
        L1 Bi-CGSTAB & 69 & 171 & 408 & 666 & 780 \\
        L2 Bi-CGSTAB & 15 & 39 & 84 & 119 & 177 \\
        L3 Bi-CGSTAB & 29 & 89 & 174 & 370 & 544 \\
        \noalign{\smallskip}\hline
    \end{tabular}}
\end{table}
\begin{table}[htbp]
    \centering
    \caption{Tolerance study for CSLP approximation for $kh=0.3125$.}
    \label{tab:tolerance-study-3125}
    \scalebox{0.8}{\begin{tabular}{cccccc}
        \hline\noalign{\smallskip}
         & $10^{-1}$ & $10^{-2}$ & $10^{-4}$ & $10^{-6}$ & $10^{-8}$ \\ \boldline
        Outer FGMRES & 6 & 5 & 5 & 5 & 5 \\
        Coarsest FGMRES & 8 & 8 & 8 & 8 & 8 \\
        L1 Bi-CGSTAB & 22 & 338 & 1149 & 1771 & 2344 \\
        L2 Bi-CGSTAB & 10 & 44 & 131 & 251 & 300 \\
        L3 Bi-CGSTAB & 34 & 48 & 114 & 198 & 272 \\
        \noalign{\smallskip}\hline
    \end{tabular}}
\end{table}

From Table \ref{tab:tolerance-study-625} and Table \ref{tab:tolerance-study-3125}, it is observed that whether the coarsest grid system remains indefinite or becomes negative definite, setting a tolerance stricter than $10^{-2}$ for the approximate inverse of CSLP does not necessarily result in a further reduction in the number of outer iterations. For cases with tolerances of $10^{-1}$ and $10^{-2}$, while the number of outer iterations is reduced by $1-2$ with a tolerance of $10^{-2}$, achieving this tolerance requires several times more iterations, particularly on the first and second levels, where computations are expensive. We choose to set the tolerance for solving the CSLP to $10^{-1}$, striking a balance between achieving sufficient accuracy in the solution and minimizing the overall computational cost. Given that the scaling behavior of CSLP is well-established in the literature \cite{erlangga2004class,erlangga2005robust}, we limit our analysis to a single configuration and wavenumber, as this adequately demonstrates the effectiveness of our chosen tolerance.

For a tolerance on the order of $10^{-1}$, where the required number of iterations is not substantial, we choose to use the more stable GMRES solver for better approximations to the inverse of CSLP. Furthermore, according to \cite{dwarka2020scalablemultilevel}, while setting the tolerance at $10^{-1}$, we limit the maximum number of iterations to $6 (N^{l})^{\frac{1}{4}}$, where $N^{l}$ denotes the size of the problem on level $l$. This allows the iterations to cease once the maximum number of iterations or the tolerance level is reached. 

\subsection{On wavenumber independent convergence}
To achieve a robust multilevel deflation method, we expand our investigation of convergence to deeper levels, larger wavenumbers, and more complex model problems.
\subsubsection{For constant wavenumber problem}
We employ the V-cycle multilevel deflation method with coarsening to different levels to solve the constant wavenumber model problem with increasing wavenumber $k$. As mentioned, GMRES instead of Bi-CGSTAB is used to solve the CSLP with complex shift $\beta_2=0.1$. The relative tolerance for solving CSLP is set to $10^{-1}$. For the coarsest problem, one iteration is performed if it is negative definite; otherwise, CSLP-preconditioned GMRES iterations are employed to reduce the relative residual to $10^{-1}$. The tolerance for the outer FGMRES iterations is $10^{-6}$. We consider the scenario with $kh=0.3125$, indicating that from the fourth level onward, the linear system becomes negative definite.

\begin{table}[htbp]
\centering
    \caption{The number of outer iterations required for the constant wavenumber model problems with increasing wavenumber $k$ by the multilevel deflation combined with CSLP with complex shift $\beta_2 = 0.1$}
    \label{tab:scalable_convergence_MP2b_beta2m01}
    \scalebox{0.8}{\begin{tabular}{cccc}
    \hline\noalign{\smallskip}
    \multirow{2}{*}{$k$} & \multicolumn{3}{c}{Multilevel Deflation} \\ 
    \noalign{\smallskip}\cline{2-4}\noalign{\smallskip}
                                  & Three-level  & Four-level  & Five-level  \\ 
    \boldline
    100                           & 6            & 9           & 8           \\
    200                           & 6            & 13          & 12          \\
    400                           & 7            & 20          & 20          \\
    800                           & 7            & 37          & 37          \\ \noalign{\smallskip}\hline
    \end{tabular}}
\end{table}

From Table \ref{tab:scalable_convergence_MP2b_beta2m01}, we observe that for the multilevel deflation method, if the coarsest grid system remains indefinite, it exhibits convergence behavior that is close to wavenumber independent, corresponding to the three-level deflation method in the table. However, if the coarsest grid system becomes negative definite, as shown by the four-level deflation method in the table, convergence results can still be achieved, but the number of outer iterations starts to increase with the wavenumber. We also find that continuing to deeper levels, as demonstrated by the five-level deflation method in the table, does not lead to an increase in the number of outer iterations compared to the four-level deflation. \hl{While theoretically we could continue to deeper levels until the coarsest problem becomes small enough for direct solving, this approach is less favorable in massively parallel computing environments due to the increased communication costs and potential load imbalance.}

\begin{table}[htbp]
\centering
    \caption{The number of outer iterations required for the constant wavenumber model problems with increasing wavenumber $k$ by the multilevel deflation combined with CSLP with complex shift $\beta_2 = {k}^{-1}$}
    \label{tab:scalable_convergence_MP2b_beta2m11byk}
    \scalebox{0.8}{\begin{tabular}{cccc}
    \hline\noalign{\smallskip}
    \multirow{2}{*}{$k$} & \multicolumn{3}{c}{Multilevel Deflation} \\ \noalign{\smallskip}\cline{2-4}\noalign{\smallskip}
                                  & Three-level  & Four-level  & Five-level  \\ \boldline
    100                           & 6            & 6           & 6           \\
    200                           & 6            & 7           & 7          \\
    400                           & 6            & 8           & 8          \\
    800                           & 7            & 9           & 9          \\ \noalign{\smallskip}\hline
    \end{tabular}}
\end{table}

As mentioned above, the Krylov subspace iterations for the CSLP allow the benefits of using a small shift, resulting in a preconditioner similar to the original Helmholtz operator that retains the ability to shift indefiniteness at certain levels. Similar to \cite{dwarka2020scalablemultilevel}, one can use the inverse of the wavenumber $k$ as the shift ($\beta_2 = k^{-1}$). As observed in Section \ref{sec:tol_for_cslp}, having a tolerance of $10^{-1}$ to approximate the inverse of CSLP leads to an increase in the number of outer iterations. For a small complex shift, the residual cannot be reduced to $10^{-1}$ within the maximum number of iterations given. Using the more stable GMRES often provides a relatively accurate approximation compared to the Bi-CGSTAB method. This is one of the reasons why GMRES is employed for approximating the inverse of CSLP on the coarse grid levels in this study.

As shown in Table \ref{tab:scalable_convergence_MP2b_beta2m11byk}, with complex shift $\beta_2 = {k}^{-1}$, the close-to wavenumber independent convergence is obtained even for the multilevel deflation methods where the coarsest grid problems become negative definite. Hereafter, we denote this configuration, which is mainly a parallel matrix-free implementation of the V-cycle multilevel deflation method proposed in \cite{dwarka2020scalablemultilevel}, as \textbf{MADP-v1} (Matrix-free multilevel Adapted Deflation Preconditioning).

\subsubsection{For non-constant wavenumber problem}
In this section, we apply MADP-v1 to non-constant wavenumber problems. Note that for the model problems described in Section \ref{sec:MPs}, due to the use of a computational domain based on actual physical dimensions rather than being scaled to a unit length, we use a complex shift $\beta_2 =(k_{\text{dim}})_{\text{max}}^{-1}$, where $k_{\text{dim}}$ is the so-called dimensionless wavenumber, defined as
\[k_{\text{dim}} = \sqrt{\left(\frac{2 \pi f}{c}\right)^2 L_x L_y}, \]
with $L_x$ and $L_y$ denoting the lengths of the computational domain in the $x$ and $y$ directions, respectively. 

In Table \ref{tab:scalable_convergence_MP3_beta2m11byk_Vcycle}, we give the results for the Wedge problem with $kh=0.349$, indicating that the linear systems become negative definite from the fourth-level coarse grid onward. We find that the latter case requires more outer iterations and significantly more CPU time. Upon further observation of the solving process, it is observed that coarsening to negative definite levels requires a higher number of GMRES iterations to approximate the CSLP compared to the scenario of coarsening to indefinite levels. In cases where the coarsening remains on indefinite levels, the tolerance of $10^{-1}$ is achieved within the maximum number of iterations. However, in cases where the coarsening goes to negative definite levels, the number of iterations reaches the maximum specified value without achieving the same tolerance. For example, consider the Wedge problem with a grid size of $1153 \times 1921$ and $f=\SI{160}{\hertz}$. On the first and second levels, the four-level deflation requires $232$ and $164$ GMRES iterations to approximate the CSLP per outer iteration, respectively, whereas the three-level deflation only requires $73$ and $49$ GMRES iterations.

Table \ref{tab:scalable_convergence_MP4_beta2m11byk_Vcycle} reports the number of iterations required and the time elapsed for the Marmousi problem with $kh=0.54$. In this scenario, the linear systems become negative definite from the third-level grid onward. Despite being coarsened to deeper negative definite levels, the number of outer iterations remains constant and the computational time is comparable. However, one can find that, for such a highly heterogeneous model problem, the number of outer iterations starts increasing with the frequency. This is consistent with the results in \cite{dwarka2020scalablemultilevel}.

\begin{table}[htbp]
\centering
    \caption{Number of iterations required and time elapsed for the Wedge problem with $kh=0.349$ for the largest wavenumber $k$. In parentheses are the number of iterations to solve the coarsest grid system. \hl{L1 and L2 represent the number of GMRES iterations required to approximate the inverse of CSLP on the first and second level, respectively.}}
    \label{tab:scalable_convergence_MP3_beta2m11byk_Vcycle}
    \scalebox{0.75}{
\begin{tabular}{ccccccccccc}
\hline\noalign{\smallskip}
    &       & \multicolumn{4}{c}{Three-level MADP-v1}   &  & \multicolumn{4}{c}{Four-level MADP-v1}    \\ 
\noalign{\smallskip}\cline{3-6} \cline{8-11}\noalign{\smallskip}
\multirow{2}{*}{$f$ ($\SI{}{\hertz}$)} & \multirow{2}{*}{Grid size} & Outer \#iter & \multicolumn{2}{c}{\#iter for CSLP }& CPU  &  & Outer \#iter & \multicolumn{2}{c}{\#iter for CSLP } & CPU  \\ 
\noalign{\smallskip}\cline{4-5} \cline{9-10}\noalign{\smallskip}
    &       & (Coarsest)   & L1 & L2        & time (s) &  & (Coarsest)   & L1 & L2       & time (s) \\ 
\boldline
20     & 145$\times$241    & 7 (2)  & 20 & 51  & 3.78     &  & 8 (1)  & 48 & 59  & 7.00      \\
40     & 289$\times$481    & 7 (3)  & 25 & 59  & 20.14    &  & 9 (1)  & 116& 83  & 103.31    \\
80     & 577$\times$961    & 8 (4)  & 38 & 61  & 195.14   &  & 11 (1) & 164& 116 & 907.00    \\
160    & 1153$\times$1921  & 8 (4)  & 73 & 49  & 1060.50  &  & 13 (1) & 232& 164 & 5101.73   \\ \noalign{\smallskip}\hline
\end{tabular}}
\end{table}
\begin{table}[htbp]
\centering
    \caption{Number of iterations required and time elapsed for the Marmousi problem with $kh=0.54$ for the largest wavenumber $k$.}
    \label{tab:scalable_convergence_MP4_beta2m11byk_Vcycle}
    \scalebox{0.8}{
\begin{tabular}{ccccccccccc}
\hline\noalign{\smallskip}
    &       & \multicolumn{4}{c}{Three-level MADP-v1}   &  & \multicolumn{4}{c}{Five-level MADP-v1}    \\ \noalign{\smallskip}\cline{3-6} \cline{8-11}\noalign{\smallskip} 
\multirow{2}{*}{$f$ ($\SI{}{\hertz}$)} & \multirow{2}{*}{Grid size} & Outer & \multicolumn{2}{c}{\#iter for CSLP}& CPU &  & Outer & \multicolumn{2}{c}{\#iter for CSLP} & CPU  \\ \noalign{\smallskip}\cline{4-5} \cline{9-10}\noalign{\smallskip} 
    &            & \#iter   & L1  & L2  & time (s)  &  & \#iter   & L1  & L2    & time (s)  \\ \boldline
10  & 737$\times$241  & 12  & 124 & 88  & 133.91    &  & 12       & 124 & 88    & 165.43   \\
20  & 1473$\times$481 & 16  & 175 & 124 & 1560.03   &  & 16       & 175 & 124   & 1743.97  \\
40  & 2945$\times$961 & 23  & 247 & 175 & 21557.94  &  & 23       & 247 & 175   & 21233.10 \\ \noalign{\smallskip}\hline
\end{tabular}}
\end{table}

In summary, the variant MADP-v1, based on the configuration proposed in \cite{dwarka2020scalablemultilevel}, utilizes a V-cycle type and allows the combination of CSLP with a smaller complex shift $\beta_2 =(k_{\text{dim}})_{\text{max}}^{-1}$. For constant wavenumber problems,  MADP-v1 achieves near wavenumber-independent convergence. However, for non-constant wavenumber model problems, as the wavenumber increases, the number of outer iterations will increase gradually, and the cost of using Krylov subspace methods to solve CSLP will become more noticeable.

For practical applications, where the wavenumber is usually non-constant within the domain, and also for better scalability, the coarsening in this multilevel deflation method should not be limited only to indefinite levels. Therefore, we will pay more attention to the common occurrence of coarsening to negative definite levels. The case of coarsening to indefinite levels will serve as a reference for the case of coarsening to negative definite levels. We aim to achieve at least similar convergence and computational efficiency in the case of coarsening to negative definite levels.

\subsubsection{On the tolerance for coarse levels}
Sheikh et al. \cite{sheikh2013convergence} stated that using $n2$, $n3$, and $n4$ iterations on the second, third, and fourth levels, respectively, can accelerate convergence, and their results indicate that larger $n2$ leads to better convergence for larger wavenumber. Additionally, \cite{dwarka2020scalablemultilevel} demonstrates that employing a W-cycle instead of a V-cycle for constructing the multilevel hierarchy results in a reduced number of iterations across reported frequencies. In this paper, we attribute this to the accuracy of solving on the second, third, and fourth levels.

In contrast to the previous configuration of a single iteration on each coarser level, we introduced distinct tolerances for iterations on the second, third, and fourth levels, exploring their impact on outer iterations and CPU time. The numerical experiments utilized the Marmousi model problem with a grid size of $1473\times961$ and frequencies of $\SI{10}{\hertz}$ and $\SI{20}{\hertz}$, corresponding to $kh=0.27$ and $0.54$, respectively. As we mentioned, for $kh=0.27$, the linear systems of the second and third levels remain indefinite, while that of the fourth level becomes negative definite. For $kh=0.54$, the third and fourth levels become negative definite.

\begin{figure}[htbp]
    \centering
    \includegraphics[width=0.8\textwidth]{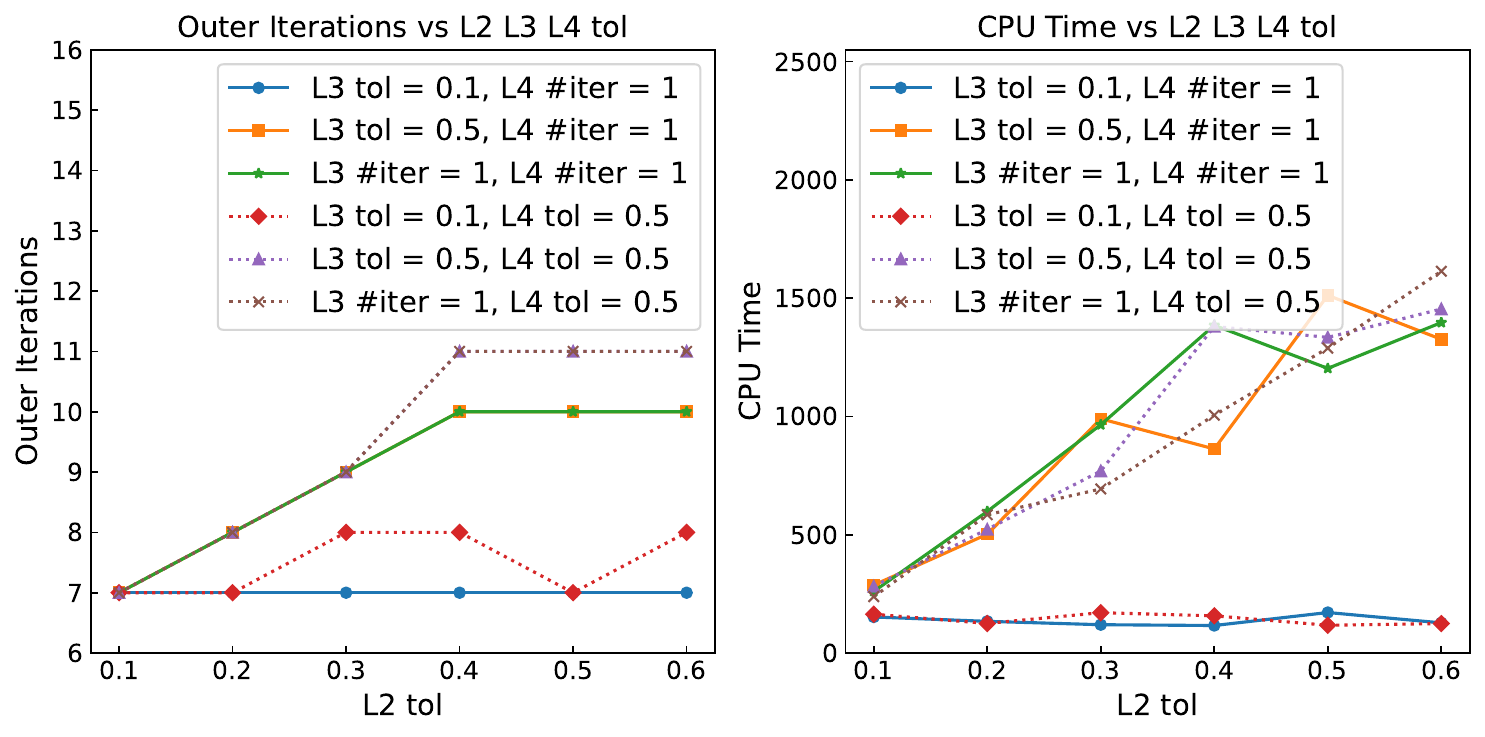}
    \caption{Outer iterations and CPU time vary from different tolerances on the second (L2), third (L3), and fourth (L4) levels. Five-level deflation for Marmousi problem, grid size $1473 \times 481$, $f=\SI{10}{\hertz}$}
    \label{fig:five_levels_nx1473f10_L2_tol_varies}
\end{figure} 

\begin{figure}[htbp]
    \centering
    \includegraphics[width=0.8\textwidth]{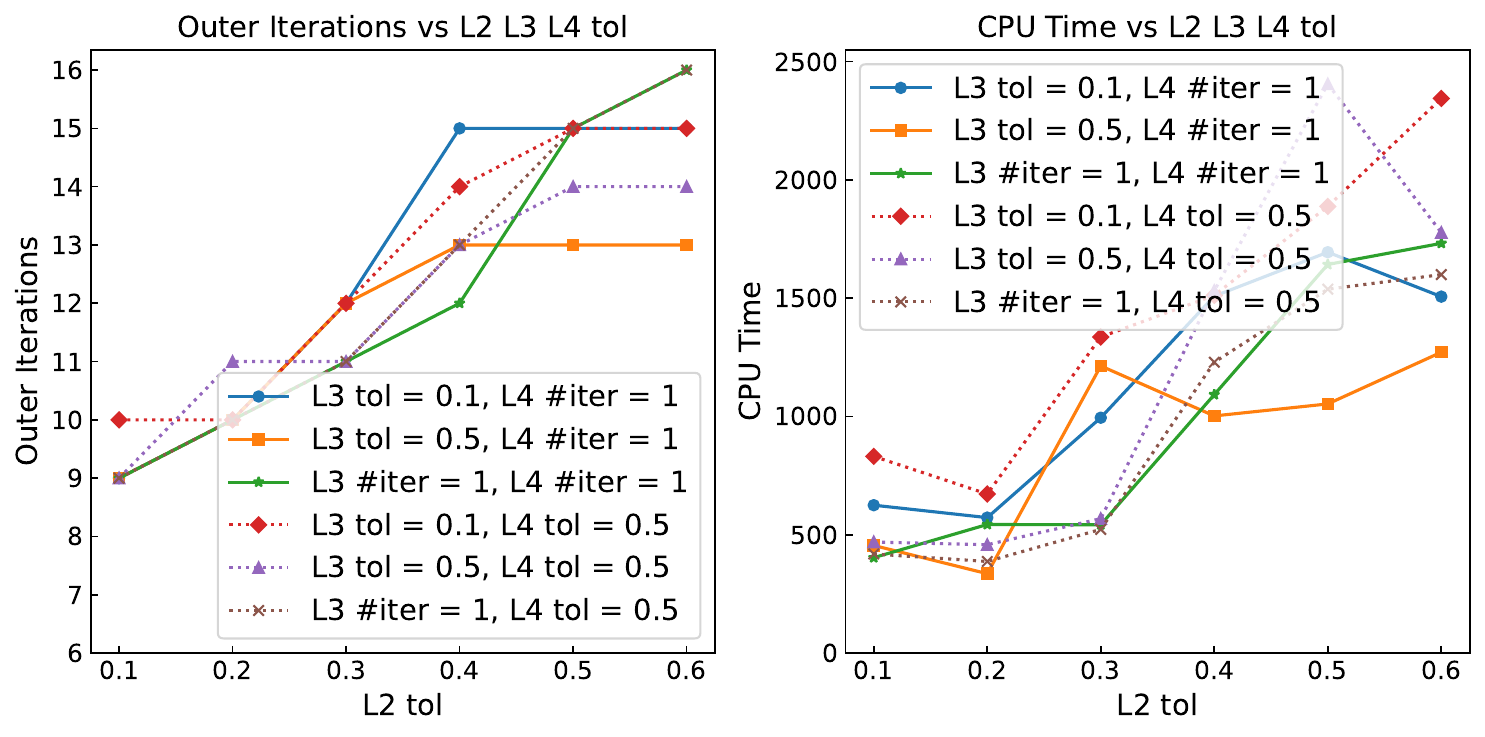}
    \caption{Outer iterations and CPU time vary from different tolerances on the second (L2), third (L3), and fourth (L4) levels. Five-level deflation for Marmousi problem, grid size $1473 \times 481$, $f=\SI{20}{\hertz}$}
    \label{fig:five_levels_nx1473f20_L2_tol_varies}
\end{figure} 

It is evident in Figure \ref{fig:five_levels_nx1473f10_L2_tol_varies} and \ref{fig:five_levels_nx1473f20_L2_tol_varies} that the number of outer iterations is correlated with the accuracy of solving the second-level grid system. Overall, a higher number of outer iterations usually corresponds to increased CPU time. However, we also observe that the minimum CPU time does not necessarily align with the minimum number of outer iterations, as depicted in Figure \ref{fig:five_levels_nx1473f20_L2_tol_varies}. This suggests that sacrificing a few extra outer iterations may result in computational time savings. A balance between the number of outer iterations and computational time needs to be identified.

For $kh=0.27$, it is best to set a tolerance of $10^{-1}$ for the second and third levels and perform one iteration for coarser levels. Conversely, for the case of $kh=0.54$, the recommended tolerances for the second and third levels are $2 \times 10^{-1}$ and $5 \times 10^{-1}$, respectively. From extensive numerical experiments across various grid sizes and multilevel deflation methods, it is observed that the optimal setting of the tolerance for solving the second, third, and fourth-level grid systems, corresponding to the minimum CPU time and the fewest outer iterations, may vary. However, on a comprehensive scale, a robust and acceptable configuration is to set a tolerance for solving the second-level grid system to $10^{-1}$, while performing only one iteration on the other coarser grid levels. Let us denote this configuration as \textbf{MADP-v2}.

We applied the MADP-v2 to solve the Wedge and Marmousi problems, respectively. The results are presented in Tables \ref{tab:scalable_convergence_MP3_beta2m11byk_L2tol} and \ref{tab:scalable_convergence_MP4_beta2m11byk_L2tol}. Compared to the corresponding results in Tables \ref{tab:scalable_convergence_MP3_beta2m11byk_Vcycle} and \ref{tab:scalable_convergence_MP4_beta2m11byk_Vcycle}, MADP-v2 results in reduced computation time and a lower number of iterations across all reported frequencies, showcasing a closer-to-wavenumber-independent behavior. In addition to reduced outer iterations, we also observe that, when setting the tolerance to $10^{-1}$ for the second-level coarse grid system, the number of iterations required to solve CSLP on the first-level grid is significantly reduced. Comparing Table \ref{tab:scalable_convergence_MP3_beta2m11byk_L2tol} with Table \ref{tab:scalable_convergence_MP3_beta2m11byk_Vcycle} (three-level deflation), in the case of coarsening to negative definite levels, MADP-v2 achieves similar convergence and computational efficiency.

\begin{table}[htbp]
\centering
    \caption{Number of outer FGMRES-iterations for the Wedge problem with $kh=0.349$ for the largest wavenumber $k$. In parentheses are the number of iterations to solve the second-level grid system. }
    \label{tab:scalable_convergence_MP3_beta2m11byk_L2tol}
    \scalebox{0.8}{
    \begin{tabular}{cccccc}
    \hline\noalign{\smallskip}
         &                  & \multicolumn{4}{c}{Four-level MADP-v2} \\ \noalign{\smallskip}\cline{3-6}\noalign{\smallskip}
    \multirow{2}{*}{$f$ ($\SI{}{\hertz}$)} & \multirow{2}{*}{Grid size} & Outer \#iter & \multicolumn{2}{c}{\#iter for CSLP}  & CPU  \\ \noalign{\smallskip}\cline{4-5}\noalign{\smallskip}
         &                  & (L2 \#iter)  & {L1}   & {L2}  &   time (s) \\ \boldline
    20   & 145$\times$241   & 6 (2)        & 6    & 59      & 4.69     \\
    40   & 289$\times$481   & 6 (2)        & 8    & 83      & 19.26    \\
    80   & 577$\times$961   & 7 (2)        & 7    & 116     & 148.05   \\
    160  & 1153$\times$1921 & 7 (3)        & 15   & 164     & 1113.86  \\ \noalign{\smallskip}\hline
    \end{tabular}}
\end{table}

\begin{table}[htbp]
\centering
    \caption{Number of outer FGMRES-iterations for the Marmousi problem with $kh=0.54$ for the largest wavenumber $k$. In parentheses are the number of iterations to solve the second-level grid system.}
    \label{tab:scalable_convergence_MP4_beta2m11byk_L2tol}
    \scalebox{0.75}{
\begin{tabular}{ccccccccccc}
\hline\noalign{\smallskip}
    &       & \multicolumn{4}{c}{Three-level MADP-v2}   &  & \multicolumn{4}{c}{Five-level MADP-v2}    \\ \noalign{\smallskip}\cline{3-6} \cline{8-11}\noalign{\smallskip} 
\multirow{2}{*}{$f$ ($\SI{}{\hertz}$)} & \multirow{2}{*}{Grid size} & Outer \#iter & \multicolumn{2}{c}{\#iter for CSLP}& CPU  &  & Outer \#iter & \multicolumn{2}{c}{\#iter for CSLP} & CPU  \\ \noalign{\smallskip}\cline{4-5} \cline{9-10}\noalign{\smallskip}
    &                  & (L2 \#iter)   & L1  & L2   & time (s)  &  & (L2 \#iter)   & L1  & L2    & time (s) \\ \boldline
10  & 737$\times$241   & 8 (3)         & 23  & 88   & 58.01     &  & 8 (3)         & 23  & 88    & 48.29   \\
20  & 1473$\times$481  & 9 (3)         & 25  & 124  & 422.33    &  & 9 (3)         & 20  & 124   & 364.21  \\
40  & 2945$\times$961  & 10 (3)        & 39  & 175  & 5267.53   &  & 10 (3)        & 39  & 175   & 4106.46 \\ \noalign{\smallskip}\hline
\end{tabular}}
\end{table}

\subsection{Combined with multigrid-based CSLP}
Further observation reveals that, in cases where coarsening reaches negative definite levels, a significant portion of the computational time is still dedicated to approximating the inverse of CSLP on the first and second levels. Moreover, these iterations on the second level typically reach the specified maximum number of iterations $6(N^l)^{\frac{1}{4}}$ rather than achieving a tolerance of $10^{-1}$. The use of GMRES iterations for solving CSLP on the first and second levels consumes a substantial amount of time, since the scale of the first- and second-level grid systems is large.

Instead of employing the GMRES or Bi-CGSTAB methods, we can utilize the multigrid method to approximate the inverse of CSLP on the first and second levels. On coarser levels, GMRES iterations are still used to approximate the inverse of CSLP. However, as known, the multigrid method requires that the complex shift should not be too small. Consequently, we cannot use $\beta_2 = (k_{\text{dim}})_{\text{max}}^{-1}$ as the complex shift for CSLP. Therefore, in the multilevel deflation methods combined with the multigrid-based CSLP, a complex shift of $\beta_2=0.5$ will be consistently utilized. (Additional numerical experiments have demonstrated that $\beta_2=0.5$ is a superior choice among other smaller complex shifts.)   

Except for the choice of the complex shift for CSLP and the method used to solve CSLP on the first- and second-level coarse grid systems, the remaining settings are mostly inherited from MADP-v2. Specifically, a tolerance of $10^{-1}$ is set for solving the second-level grid system, and only one iteration is performed on other coarser levels. This modified configuration is denoted as \textbf{MADP-v3}. As it combines with multigrid-based CSLP on the first and second levels, this variant can be considered as an extension of the two-level deflation method proposed in \cite{chen2023matrixfree2ldef}.

The number of iterations and computation time required for solving the Wedge and Marmousi problems using MADP-v3 are presented in Tables \ref{tab:scalable_convergence_MP3_L12CSLPbeta2m05_L2tol} and \ref{tab:scalable_convergence_MP4_L12CSLPbeta2m05_L2tol}, respectively. We observe that compared to MADP-v2 (as shown in Tables \ref{tab:scalable_convergence_MP3_beta2m11byk_L2tol} and \ref{tab:scalable_convergence_MP4_beta2m11byk_L2tol}), while the number of outer iterations has increased, it exhibits nearly wavenumber-independent convergence, with computation time three times faster. Moreover, compared to the two-level deflation method proposed in  \cite{chen2023matrixfree2ldef}, the current multilevel deflation method ensures a similar number of outer iterations while significantly reducing computation time.

\begin{table}[htbp]
\centering
    \caption{Number of outer FGMRES-iterations for the Wedge problem with $kh=0.349$ for the largest wavenumber $k$. In parentheses are the number of iterations to solve the second-level grid system. }
    \label{tab:scalable_convergence_MP3_L12CSLPbeta2m05_L2tol}
    \scalebox{0.8}{
    \begin{tabular}{ccccccc}
    \hline\noalign{\smallskip}
        &       & \multicolumn{2}{c}{Four-level MADP-v3}   &  & \multicolumn{2}{c}{Five-level MADP-v3}    \\ \noalign{\smallskip}\cline{3-4} \cline{6-7}\noalign{\smallskip}
    \multirow{2}{*}{$f$ ($\SI{}{\hertz}$)} & \multirow{2}{*}{Grid size} & Outer \#iter & CPU  &  & Outer \#iter  & CPU  \\
        &                  & (L2 \#iter)   & time (s)  &  & (L2 \#iter)   & time (s) \\ \boldline
    20   & 145$\times$241   & 10 (6)       & 1.73      &  & 10 (6)        & 1.83 \\
    40   & 289$\times$481   & 10 (10)      & 8.08      &  & 10 (10)       & 8.87 \\
    80   & 577$\times$961   & 10 (17)      & 48.05     &  & 10 (18)       & 64.54 \\
    160  & 1153$\times$1921 & 11 (34)      & 356.76    &  & 11 (34)       & 367.53 \\
    320  & 2305$\times$3841 & 11 (66)      & 3458.14   &  & 11 (64)       & 3065.03 \\
    \noalign{\smallskip}\hline
    \end{tabular}}
\end{table}
\begin{table}[htbp]
\centering
    \caption{Number of outer FGMRES-iterations for the Marmousi problem with $kh=0.54$ for the largest wavenumber $k$. In parentheses are the number of iterations to solve the second-level grid system.}
    \label{tab:scalable_convergence_MP4_L12CSLPbeta2m05_L2tol}
    \scalebox{0.8}{
\begin{tabular}{cccccccccc}
\hline\noalign{\smallskip}
    &       & \multicolumn{2}{c}{Two-level Deflation\cite{chen2023matrixfree2ldef}}& &\multicolumn{2}{c}{Three-level MADP-v3}   &  & \multicolumn{2}{c}{Five-level MADP-v3}    \\ \noalign{\smallskip}\cline{3-4} \cline{6-7} \cline{9-10}\noalign{\smallskip}
\multirow{2}{*}{$f$ ($\SI{}{\hertz}$)} & \multirow{2}{*}{Grid size} &Outer & CPU & & Outer \#iter & CPU &  & Outer \#iter  & CPU\\
    &              & \#iter & time (s) & & (L2 \#iter)   &  time (s)  &  & (L2 \#iter)   & time (s)  \\ \boldline
10  & 737$\times$241   & 11   & 23.15      &  & 11 (17)       & 16.53     &  & 11 (13)       & 18.57   \\
20  & 1473$\times$481  & 11   & 224.21     &  & 11 (30)       & 110.79    &  & 11 (24)       & 108.03  \\
40  & 2945$\times$961  & 12   & 4354.83    &  & 13 (69)       & 1220.61   &  & 13 (50)       & 1084.42 \\ \noalign{\smallskip}\hline
\end{tabular}}
\end{table}

\begin{figure}[htbp]
    \centering
    \includegraphics[width=0.8\textwidth]{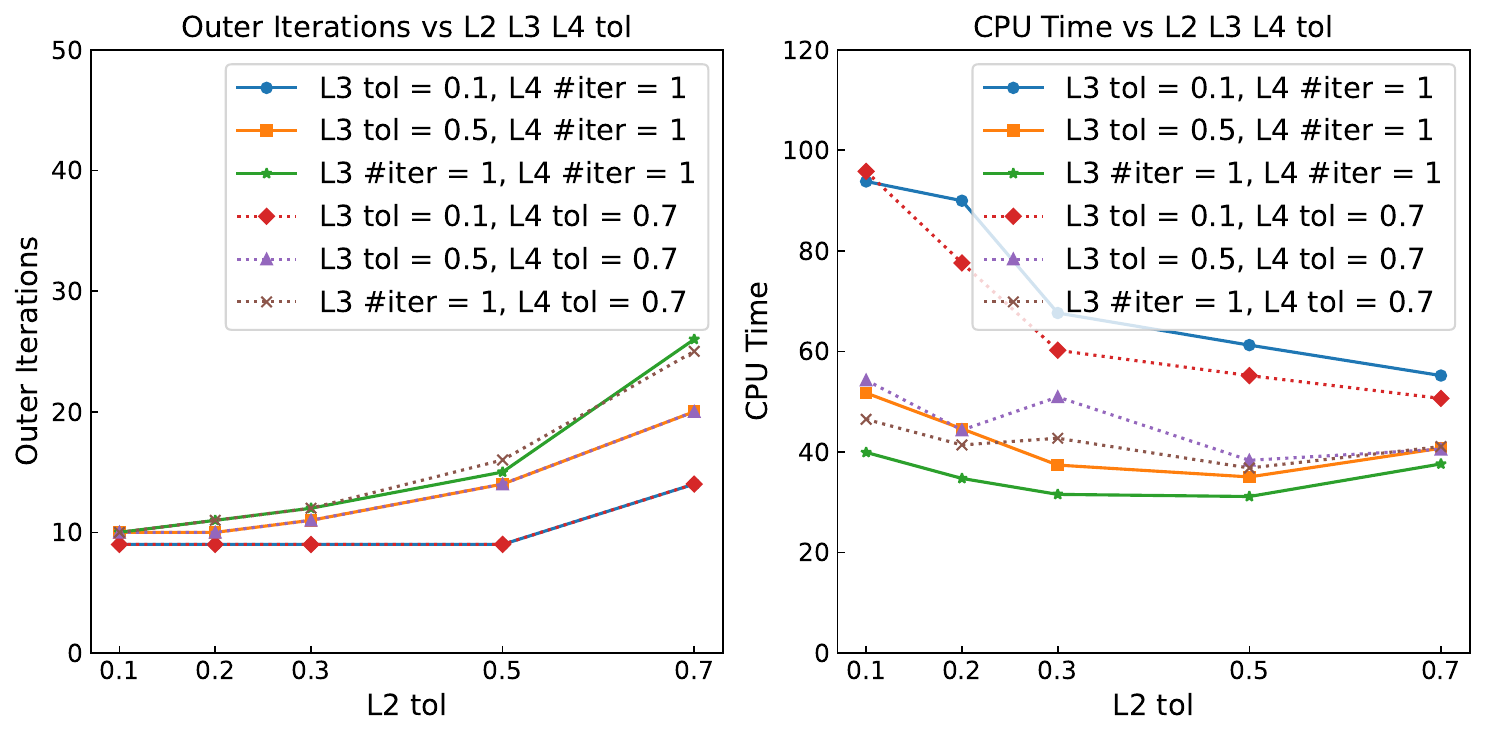}
    \caption{Outer iterations and CPU time vary from different tolerances on the second (L2), third (L3), and fourth (L4) levels. Five-level deflation combined with multigrid-based CSLP on first and second levels. Marmousi problem, grid size $1473 \times 481$, $f=\SI{10}{\hertz}$}
    \label{fig:five_levels_L1L2MG05_nx1473f10_L2_tol_varies}
\end{figure}

\begin{figure}[htbp]
    \centering
    \includegraphics[width=0.8\textwidth]{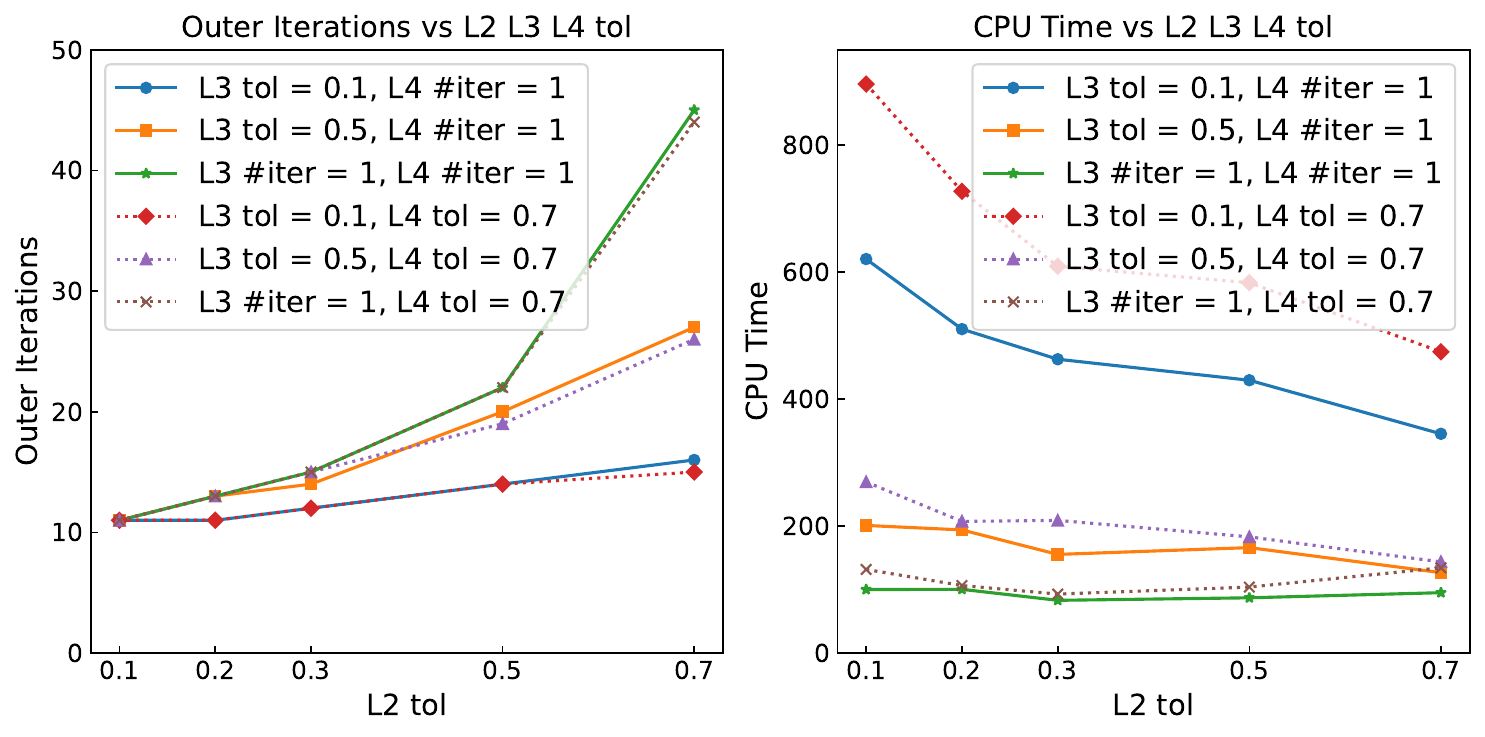}
    \caption{Outer iterations and CPU time vary from different tolerances (tol) on the second (L2), third (L3), and fourth (L4) levels. Five-level deflation combined with multigrid-based CSLP on first and second levels. Marmousi problem, grid size $1473 \times 481$, $f=\SI{20}{\hertz}$}
    \label{fig:five_levels_L1L2MG05_nx1473f20_L2_tol_varies}
\end{figure}

Similarly to the last section, the optimal tolerance setting is studied for the second, third and fourth levels, as shown in Figures \ref{fig:five_levels_L1L2MG05_nx1473f10_L2_tol_varies} and \ref{fig:five_levels_L1L2MG05_nx1473f20_L2_tol_varies}. From the figures, it can be observed that performing only one iteration on the fourth level (L4) or setting a tolerance of $0.5$ has little impact on the outer iterations but introduces additional computational costs. For this reason, we can keep one iteration on the fourth level. In comparison to performing only one iteration on the third level (L3), setting a smaller tolerance on L3 helps slow down the increase in the number of outer iterations but leads to more computational costs. From the perspective of computation time, performing only one iteration on L3 remains the optimal choice. With the incorporation of multigrid-based CSLP, the number of outer iterations increases as the tolerance of the second level (L2) increases, while the computation time shows a trend of decreasing first and then increasing. If one wants to minimize the computation time, choosing the tolerance of the second level as $0.3$ while still performing only one iteration on other coarser levels can be the optimal option. Let us denote this configuration as \textbf{MADP}. 

Tables \ref{tab:scalable_convergence_MP3_L12CSLPbeta2m05_L2tol03} and \ref{tab:scalable_convergence_MP4_L12CSLPbeta2m05_L2tol03} present the number of iterations required and computation time to solve the Wedge and Marmousi problems using MADP. Compared with Tables \ref{tab:scalable_convergence_MP3_L12CSLPbeta2m05_L2tol} and \ref{tab:scalable_convergence_MP4_L12CSLPbeta2m05_L2tol}, it can be seen that although a few extra outer iterations are consumed, reduced computation time is obtained.
\begin{table}[htbp]
\centering
    \caption{Number of outer FGMRES-iterations for the Wedge problem with $kh=0.349$ for the largest wavenumber $k$. In parentheses are the number of iterations to solve the second-level grid system. }
    \label{tab:scalable_convergence_MP3_L12CSLPbeta2m05_L2tol03}
    \scalebox{0.8}{
    \begin{tabular}{ccccccc}
    \hline\noalign{\smallskip}
        &       & \multicolumn{2}{c}{Four-level MADP}   &  & \multicolumn{2}{c}{Five-level MADP}    \\ \noalign{\smallskip}\cline{3-4} \cline{6-7}\noalign{\smallskip}
    \multirow{2}{*}{$f$ ($\SI{}{\hertz}$)} & \multirow{2}{*}{Grid size} & Outer \#iter & CPU  &  & Outer \#iter  & CPU \\
        &                  & (L2 \#iter)   & time (s)  &  & (L2 \#iter)   & time (s) \\ \boldline
    20   & 145$\times$241   & 11 (4)       & 1.51      &  & 11 (4)        & 1.55    \\
    40   & 289$\times$481   & 12 (6)       & 6.34      &  & 12 (6)        & 6.42    \\
    80   & 577$\times$961   & 13 (10)      & 39.62     &  & 13 (10)       & 42.21    \\
    160  & 1153$\times$1921 & 15 (17)      & 410.27    &  & 14 (18)       & 400.64   \\
    320  & 2305$\times$3841 & 16 (33)      & 2748.70   &  & 16 (32)       & 2762.84 \\\noalign{\smallskip}\hline
    \end{tabular}}
\end{table}

\begin{table}[htbp]
\centering
    \caption{Number of outer FGMRES-iterations for the Marmousi problem with $kh=0.54$ for the largest wavenumber $k$. In parentheses are the number of iterations to solve the second-level grid system.}
    \label{tab:scalable_convergence_MP4_L12CSLPbeta2m05_L2tol03}
    \scalebox{0.8}{
\begin{tabular}{ccccccc}
\hline\noalign{\smallskip}
    &       & \multicolumn{2}{c}{Three-level MADP}   &  & \multicolumn{2}{c}{Five-level MADP}    \\ \noalign{\smallskip}\cline{3-4} \cline{6-7}\noalign{\smallskip} 
\multirow{2}{*}{$f$ ($\SI{}{\hertz}$)} & \multirow{2}{*}{Grid size} & Outer \#iter & CPU &  & Outer \#iter  & CPU  \\
    &                  & (L2 \#iter)   & time (s)   &  & (L2 \#iter) & time (s) \\ \boldline
10  & 737$\times$241   & 14 (10)       & 12.42    &  & 13 (7)        & 12.67 \\
20  & 1473$\times$481  & 16 (19)       & 93.08    &  & 15 (15)       & 84.06 \\
40  & 2945$\times$961  & 19 (43)       & 929.62   &  & 18 (29)       & 816.38 \\ \noalign{\smallskip}\hline
\end{tabular}}
\end{table}

Therefore, we regard the variant \textbf{MADP}, which can balance both convergence and computational efficiency, as an optimal configuration of the present matrix-free multilevel deflation method. This variant employs a tolerance of $0.3$ for solving the second-level grid system and performs only one iteration on other coarser levels. A multigrid V-cycle is used to solve the CSLP on the first- and second-level grid systems. On coarser levels, several GMRES iterations approximate the inverse of CSLP with a tolerance of $10^{-1}$. In the subsequent sections, we will use this variant for numerical experiments.

\subsection{Complexity Analysis}
We next analyze the complexity of the present multilevel deflation method in relation to the problem size or to the frequency, equivalently. In this numerical experiment, the Wedge model problem is solved using \textbf{five-level MADP}. The grid resolution, \textit{i.e.} $kh$, is kept fixed to a specific value, while the frequency is growing from $\SI{10}{\hertz}$ to $\SI{160}{\hertz}$ for $kh=0.349$ or even $\SI{640}{\hertz}$ for $kh=0.1745$, respectively. The case of $f=\SI{640}{\hertz}$ leads to a linear system with approximately 142 million unknowns.
The number of outer iterations and the number of iterations on the second level are reported in Table \ref{tab:scalable_convergence_MP3_L12CSLPbeta2m05_L2tol03_kh017}. Similarly to the results in Table \ref{tab:scalable_convergence_MP3_L12CSLPbeta2m05_L2tol03}, the number of outer iterations is rather moderate and is found to grow slightly with respect to frequency. 

\begin{figure}[htbp]
    \centering
    \includegraphics[width=0.6\textwidth]{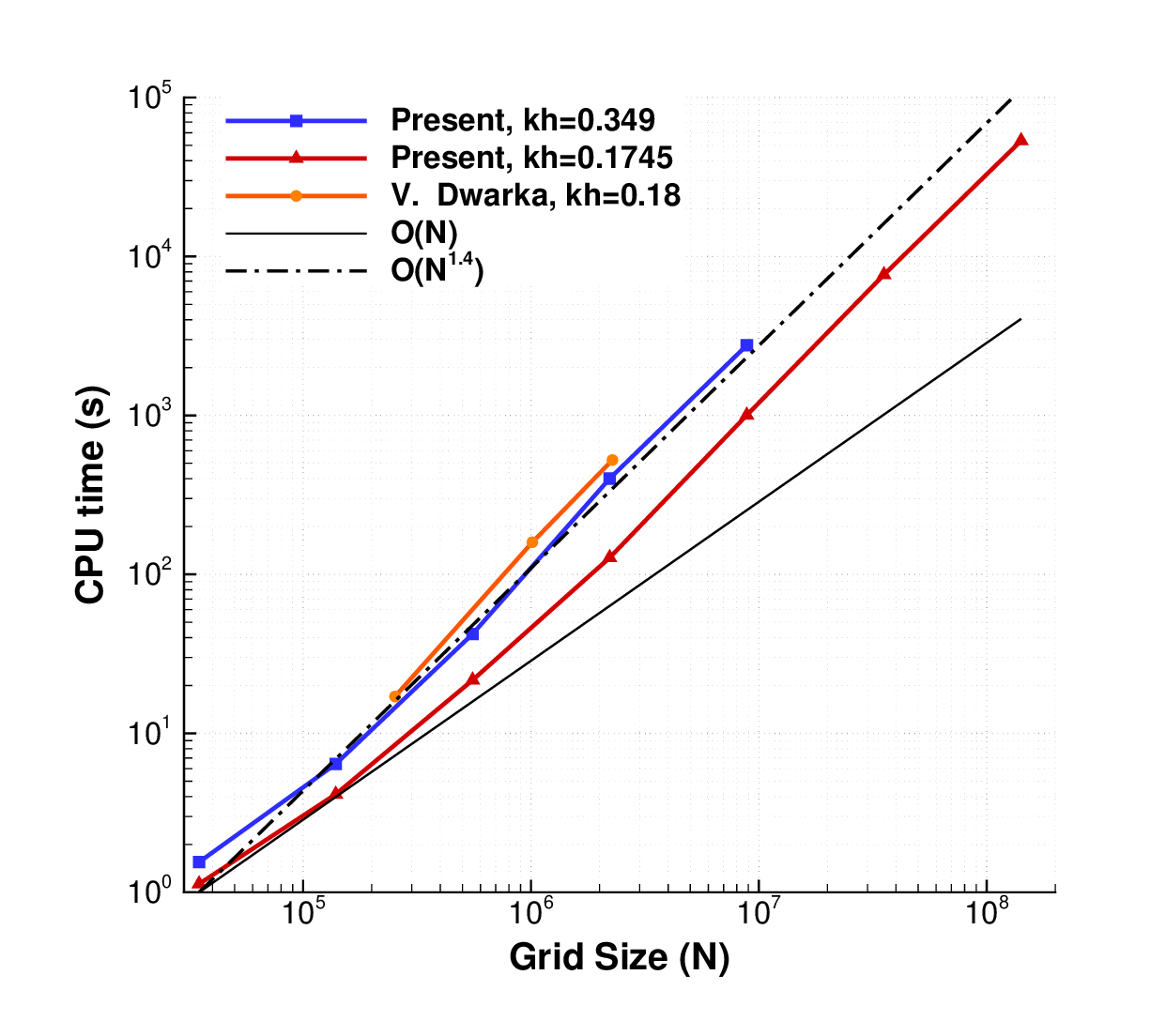}
    \caption{Evolution of computational time versus problem size. Wedge model problem. The data in orange is extracted from \cite{dwarka2020scalable}.}
    \label{fig:Wedge_complexity_analysis}
\end{figure}

Figure \ref{fig:Wedge_complexity_analysis} shows the evolution of computational time versus problem size for $kh=0.349$ and $kh=0.1745$. 
As the grid size increases, the computational time of the present matrix-free multilevel deflation method shows a similar trend to the matrix-based version proposed by \cite{dwarka2020scalablemultilevel}. However, with single-core sequential computing, the present method can handle grid sizes much larger than those achievable in \cite{dwarka2020scalablemultilevel}. If $N$ represents the total number of unknowns, it has been observed that the computational time follows a behavior of $\mathcal{O}(N)$ for small grid sizes and asymptotically approaches $\mathcal{O}(N^{1.4})$. This is comparable to the geometric two-grid preconditioner in \cite{calandra2017geometric}. The reason for this behavior is that, as the frequency increases, the number of iterations required on the second level almost increases linearly with frequency, as shown in Table \ref{tab:scalable_convergence_MP3_L12CSLPbeta2m05_L2tol03_kh017}.
\begin{table}[htbp]
      \centering
      \caption{Number of outer FGMRES-iterations for the Wedge problem with $kh=0.1745$. In parentheses are the number of iterations to solve the second-level grid system. }
      \label{tab:scalable_convergence_MP3_L12CSLPbeta2m05_L2tol03_kh017}
      \scalebox{0.8}{
      \begin{tabular}{crccc}
      \hline\noalign{\smallskip}
      \multirow{2}{*}{Grid size} & \multirow{2}{*}{\#unknowns} & \multirow{2}{*}{$f$ ($\SI{}{\hertz}$)} & Outer \#iter   & \multirow{2}{*}{CPU time (s)}\\
        &                &         & (L2 \#iter)  & \\
      \boldline
      145$\times$ 241      & 34945      & 10    & 10 (2)   &     1.13 \\
      289$\times$ 481      & 139009     & 20    & 11 (3)   &     4.14 \\
      577$\times$ 961      & 554497     & 40    & 12 (4)   &    21.64 \\
      1153$\times$ 1921    & 2214913    & 80    & 12 (7)   &   127.47  \\
      2305$\times$ 3841    & 8853505    & 160   & 13 (13)  &  1003.71  \\
      4609$\times$ 7681    & 35401729   & 320   & 14 (27)  &  7678.83  \\
      9217$\times$ 15361   & 141582337  & 640   & 17 (47)  & 53481.69 \\
      \noalign{\smallskip}\hline
      \end{tabular}}
\end{table}

\begin{figure}[htbp]
    \centering
    \includegraphics[width=0.6\textwidth]{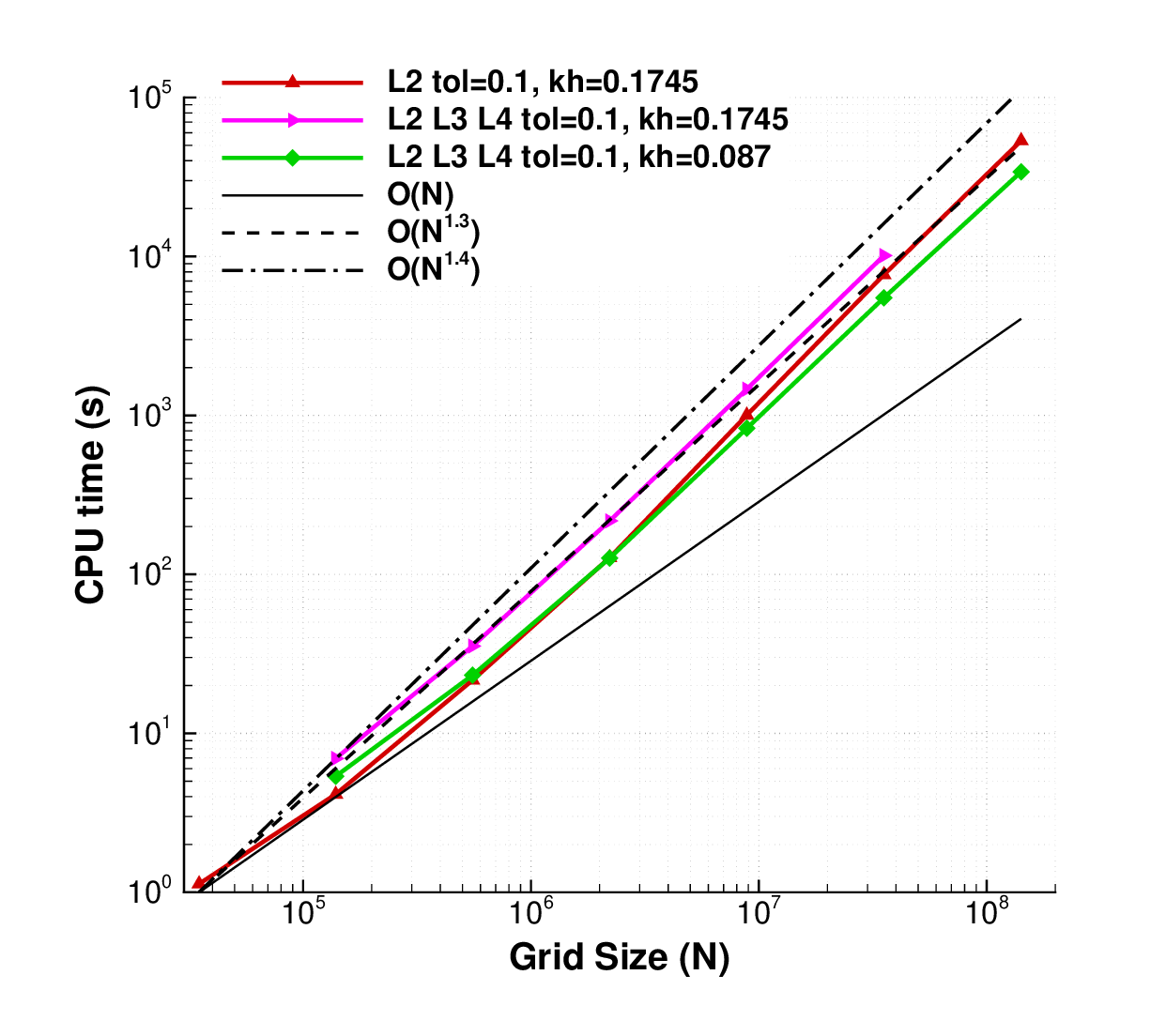}
    \caption{Evolution of computational time versus problem size. Wedge model problem.}
    \label{fig:Wedge_complexity_analysis_l234to01}
\end{figure}

One can think about continuing a similar approach, setting a tolerance of $10^{-1}$ on the third level, ensuring that the number of iterations on the second level is independent of the wavenumber, and so forth. This is feasible but only limited to indefinite levels. Setting tolerance on negative definite levels, \textit{i.e.}, performing more than one iteration, may lead to a significant increase in outer iterations and computational time, consistent with the conclusion in Section \ref{sec:coarsest_tol}. For instance, considering the Wedge model problem with $kh=0.1745$, where the fourth-level grid system remains indefinite, turning negative definite onwards the fifth level. We can extend MADP-v3 by setting a tolerance of $10^{-1}$ for iterations on the third and fourth levels instead of performing one iteration. Table \ref{tab:scalable_convergence_MP3_L12CSLPbeta2m05_L234tol01_kh017} provides the required number of outer iterations and the number of iterations on the second, third, and fourth levels. We can observe that the number of outer iterations and iterations on the second and third levels are almost independent of the wavenumber, while the number of iterations on the fourth level gradually increases with the wavenumber. However, as shown in Figure \ref{fig:Wedge_complexity_analysis_l234to01}, the computation time required is more than MADP. Together with the case of $kh=0.087$ in the figure, it can be observed that there are subtle differences in the growth trend compared to that of MADP, possibly approaching $\mathcal{O}(N^{1.3})$. While it is beneficial to set a tolerance for coarser levels for problems with smaller $kh$, whether the present multilevel deflation method can be closer to $\mathcal{O}(N)$ remains an open problem that requires further study. To complement this study, it would be interesting to perform the same complexity analysis for more levels and three-dimensional cases. These are left to a future line of research.  Additionally, in the numerical experiments on current two-dimensional problems, it does not lead to a significant reduction in computation time. Therefore, we consider MADP to still be the optimal choice.

\begin{table}[htbp]
      \centering
      \caption{Number of outer iterations and the number of iterations on the second, third, and fourth levels when a tolerance of $10^{-1}$ is set on these levels. Wedge problem with $kh=0.1745$. }
      \label{tab:scalable_convergence_MP3_L12CSLPbeta2m05_L234tol01_kh017}
      \scalebox{0.8}{
      \begin{tabular}{cccccc}
      \hline\noalign{\smallskip}
      Grid size & $f$ ($\SI{}{\hertz}$) & Outer \#iter  & L2 \#iter & L3 \#iter &L4 \#iter\\
      \boldline
      289$\times$ 481         & 20    & 10   &  3  &  2 & 16 \\
      577$\times$ 961         & 40    & 10   &  3  &  2 & 20 \\
      1153$\times$ 1921       & 80    & 10   &  3  &  2 & 30 \\
      2305$\times$ 3841       & 160   & 10   &  3  &  2 & 46 \\
      4609$\times$ 7681       & 320   & 9    &  3  &  2 & 75 \\
      \noalign{\smallskip}\hline
      \end{tabular}}
\end{table}

\section{Parallel performance}
In this section, we aim to present both weak scalability and strong scalability. Through this analysis, our goal is to offer insight into the suitability of the present multilevel deflation method for practical large-scale applications in the context of heterogeneous time-harmonic wave problems.

The \textbf{parallel six-level MADP} preconditioned FGMRES is used as the default approach in this section to solve model problems. All numerical experiments are carried out on the Linux supercomputer DelftBlue \cite{DHPC2022}, which operates on the Red Hat Enterprise Linux 8 operating system. Each compute node is furnished with two Intel Xeon E5-6248R CPUs featuring 24 cores at 3.0 GHz, 192 GB of RAM, a memory bandwidth of 132 GByte/s per socket, and a 100 Gbit/s InfiniBand card. The present solver is developed in Fortran 90 and is compiled using GNU Fortran 8.5.0 with the compiler options \verb|-O3| for optimization purposes. The Open MPI library (v4.1.1) is used for data communication.

\subsection{Weak scalability}
To assess the weak scalability of the proposed matrix-free parallel multilevel deflation preconditioning method, we keep the wavenumber or frequency unchanged and solve the model problems across varying problem sizes but maintain a fixed workload per processor. The computational times for different problem sizes and the corresponding number of processors are summarized in Tables \ref{tab:weak_scaling_MP2_L12CSLPbeta2m05_L2tol03} and \ref{tab:weak_scaling_Wedge_L12CSLPbeta2m05_L2tol03}. As the grid undergoes refinement while maintaining a constant wavenumber, the parameter $kh$ gradually decreases. In the context of deflation preconditioning, it has been documented that a smaller
$kh$ leads to a reduction in the number of outer iterations \cite{sheikh2013convergence}. As $kh$ continues to diminish, the number of outer iterations tends to stabilize. Additionally, the advantages of one or two less iterations may be counteracted by the overhead of data communication. Consequently, we observe that the computational time initially decreases due to the reduced number of outer iterations, and then remains almost constant, even as the grid size expands to tens of millions with over a thousand parallel computing cores.

\begin{table}[htbp]
      \centering
      \caption{Weak scaling for the model problem with constant wavenumber.}
      \label{tab:weak_scaling_MP2_L12CSLPbeta2m05_L2tol03}
      \scalebox{0.8}{
      \begin{tabular}{crccc}
        \hline\noalign{\smallskip}
        Grid size   & \#unknowns & np   & \#iter & CPU time (s) \\ \boldline
        \multicolumn{5}{c}{$k=400$}                  \\
         641$\times$641    & 410881     & 1    & 16     & 49.68        \\
        1281$\times$1281   & 1640961    & 4    & 13     & 21.63        \\
        2561$\times$2651   & 6558721    & 16   & 12     & 16.13        \\
        5121$\times$5121   & 26224641   & 64   & 11     & 21.66        \\
                    & &     &        &              \\
        \multicolumn{5}{c}{$k=1600$}                 \\
         2561$\times$2561  & 6558721    & 16   & 20     & 168.26       \\
         5121$\times$5121  & 26224641   & 64   & 14     & 100.84       \\
        10241$\times$10241 & 104878081  & 256  & 13     & 79.69        \\
        20481$\times$20481 & 419471361  & 1024 & 13     & 93.62        \\ \noalign{\smallskip}\hline
        \end{tabular}}
\end{table}
\begin{table}[htbp]
    \centering
    \caption{Weak scaling for the Wedge model problem with $f=\SI{320}{\hertz}$.}
    \label{tab:weak_scaling_Wedge_L12CSLPbeta2m05_L2tol03}
    \scalebox{0.8}{
        \begin{tabular}{crccc}
        \hline\noalign{\smallskip}
        Grid size  & \#unknowns& np  & \#iter & CPU time (s) \\ \boldline
        2305$\times$3841  & 8853505       & 48  & 16     & 69.75        \\
        4609$\times$7681  & 35401729      & 192 & 14     & 53.20        \\
        9217$\times$15361 & 141582337     & 768 & 14     & 67.03        \\ \noalign{\smallskip}\hline
        \end{tabular}
    }
\end{table}

This behavior is highly commendable, as it allows for the efficient resolution of large linear systems within a reasonable computational timeframe on a parallel distributed memory machine. It is important to emphasize that the advantages of the suggested approach should be considered in the context of minimizing pollution error by grid refinement for real-world application of Helmholtz problems.

\subsection{Strong scalability}
We are also interested in the strong scalability properties of the present parallel multilevel deflation preconditioning method for the Helmholtz problems. In this section, we perform numerical experiments on the nonconstant-wavenumber model problems with fixed problem sizes while varying the number of processors. First of all, the numerical experiments show that the number of iterations required is found to be independent of the number of processors used for parallel computing, which is a favorable property of our multilevel deflation method. 
\begin{figure}[htbp]
     \centering
     \begin{subfigure}[ht]{0.45\textwidth}
         \centering
         \includegraphics[width=\textwidth]{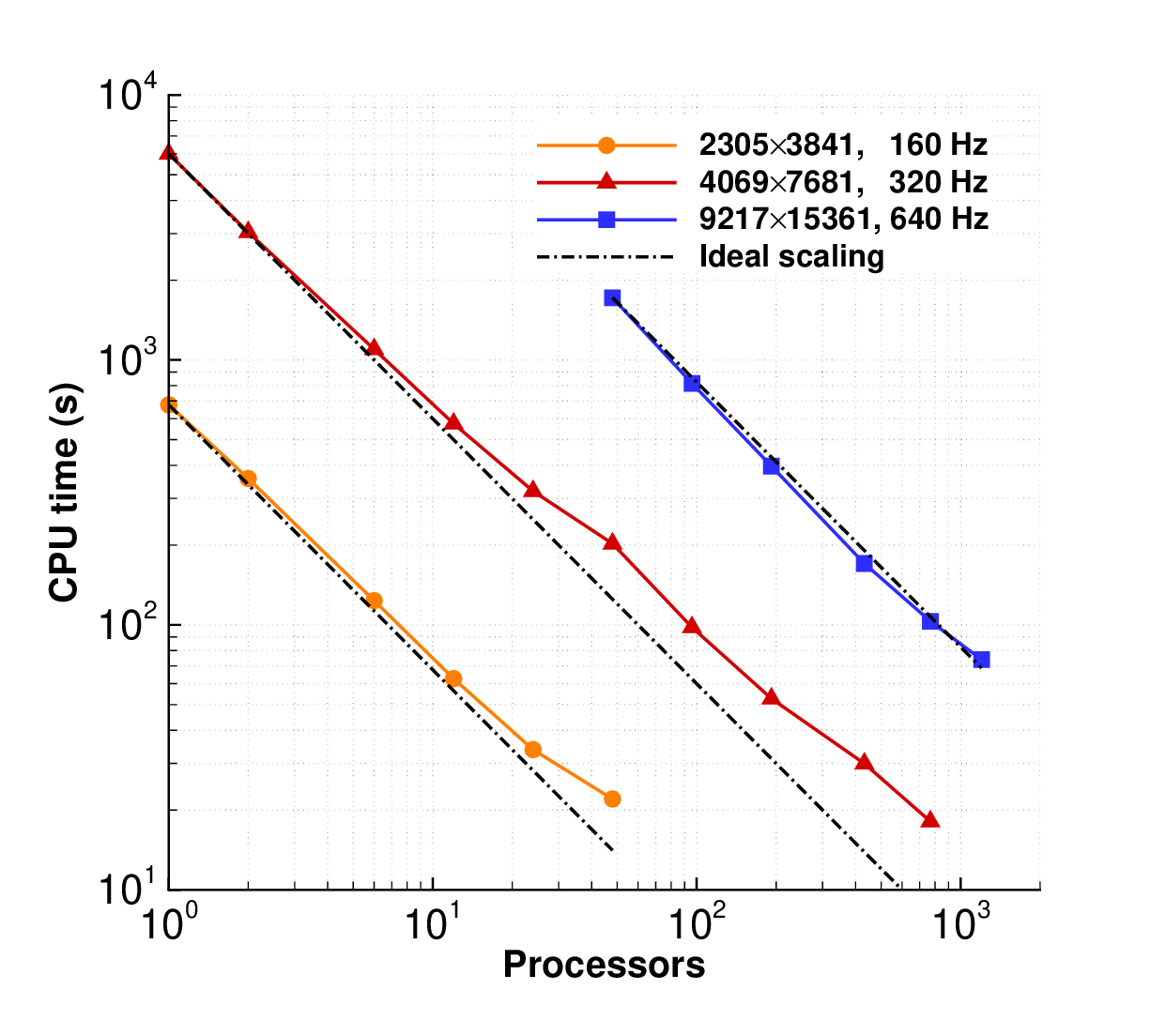}
         \caption{Wedge problem}
         \label{fig:MLDEF_L2tol03_wedge_strong}
     \end{subfigure}
     \hfill
     \begin{subfigure}[ht]{0.45\textwidth}
         \centering
         \includegraphics[width=\textwidth]{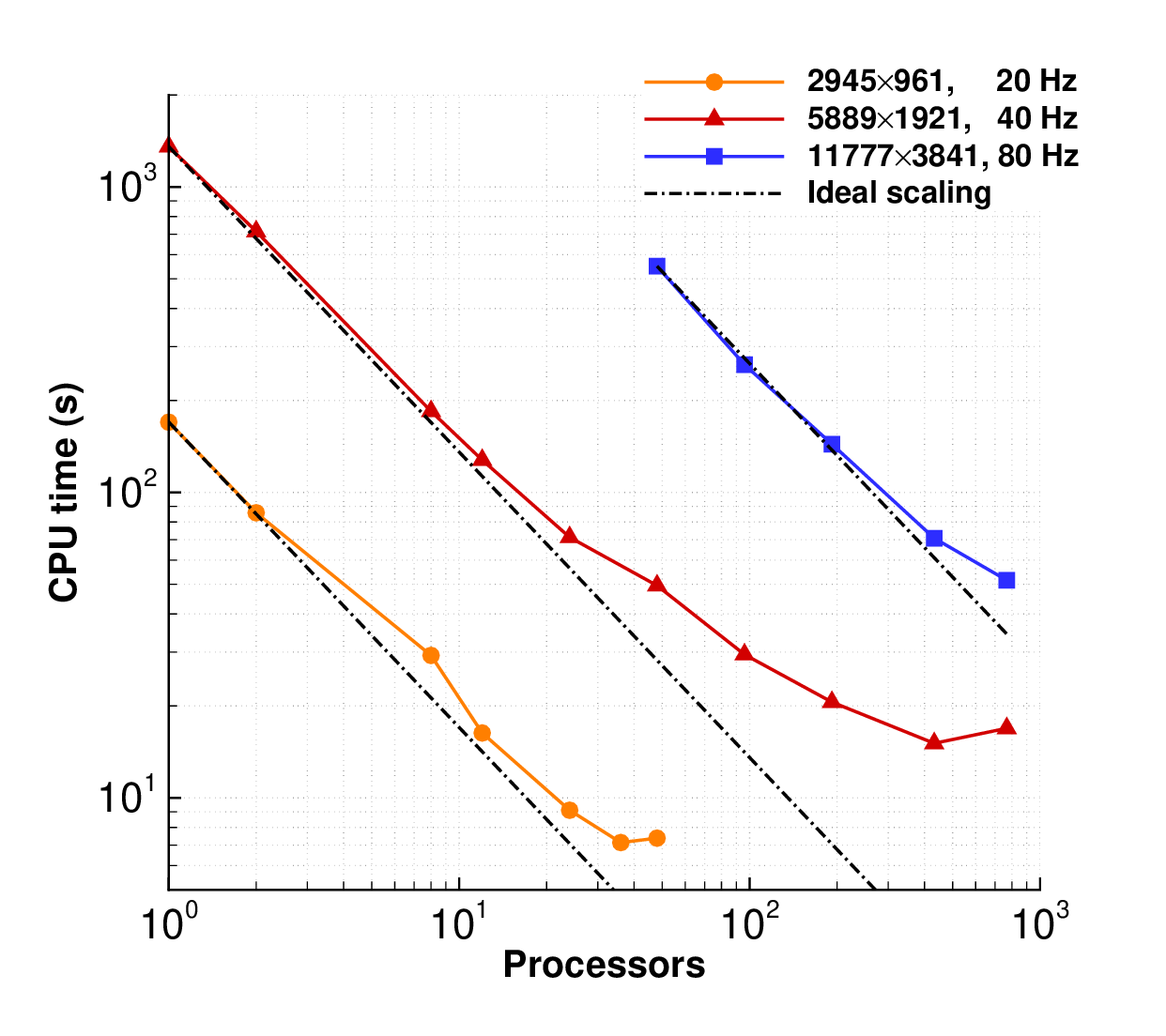}
         \caption{Marmousi problem}
         \label{MLDEF_L2tol03_marmousi_strong}
     \end{subfigure}
        \caption{Strong scaling of the parallel multilevel APD-FGMRES for the non-constant wavenumber model problems with various grid sizes and frequencies.}
    \label{fig:strong-scaling-MLDEF_L2tol03}
\end{figure}

Figure \ref{fig:strong-scaling-MLDEF_L2tol03} plots the computational time versus the number of processors. The figures show a decrease in parallel efficiency as the number of processors increases, particularly for the Marmousi model. The analysis suggests that maintaining a minimum of around one million unknowns per processor ensures a parallel efficiency of 60\% or higher. If there are fewer than 50 thousand unknowns per processor, the ratio of computation load to data communication may significantly decrease, leading to poor parallel efficiency. However, increasing the grid size for the same model problem can result in improved parallel efficiency. This is because, while the number of ghost-grid layers used for communication remains constant, the amount of data to be communicated doubles when the number of grid points doubles in each direction. Meanwhile, the total number of grid points increases fourfold, resulting in a larger ratio of computation load to data communication and better parallel efficiency. Overall, as demonstrated in solving the Wedge problem with a grid size of $9217 \times 15361$ and a frequency of $f=\SI{640}{\hertz}$, the current matrix-free multilevel deflation approach can effectively solve complex Helmholtz problems with grid sizes up to tens of millions. It demonstrates strong parallel scalability, maintaining efficiency across more than a thousand processors.

\section{Conclusion}
In this work, we present an advanced matrix-free parallel scalable multilevel deflation preconditioning method for solving the Helmholtz equation in heterogeneous time-harmonic wave problems, benchmarked on large-scale real-world models. Building on recent advancements in higher-order deflation preconditioning, our approach extends these techniques to a parallel implementation. The incorporation of the deflation technique with CSLP, along with higher-order deflation vectors and re-discretization schemes derived from the Galerkin coarsening approach, forms a comprehensive setup for matrix-free parallel implementation. The proposed re-discretized finite-difference schemes at each coarse level contribute to a convergence behavior similar to that of the matrix-based deflation method.

We have explored different configurations of the multilevel deflation method, conducting research and comparing various variants. We note that the performance of different cycle types of the present multilevel deflation method is impacted by whether the coarsest-level system is negative definite. We suggest that for the multilevel deflation method coarsening to negative definite levels, ensuring a certain accuracy in iterations on the second-level grid is crucial to maintain a consistent number of outer iterations. Based on the complexities revealed during our study, we propose a robust and efficient variant, MADP. This variant employs the following settings: a tolerance of $0.3$ for solving the second-level grid system, with only one iteration performed on other coarser levels; a multigrid V-cycle to solve CSLP on the first- and second-level grid systems; several GMRES iterations to approximate the inverse of CSLP on coarser levels, using a tolerance of $10^{-1}$. It addresses the challenges posed by negative definite coarsest-level systems and does not lead to worse complexities.

Our numerical experiments illustrate the effectiveness of the matrix-free parallel multilevel deflation preconditioner, demonstrating convergence properties that are nearly independent of the wavenumber. The reduction in memory consumption achieved through matrix-free implementation, along with satisfactory weak and strong parallel scalability, emphasizes the practical applicability of our approach for large-scale real-world applications in wave propagation.

\clearpage
\appendix

\section{Re-discretization scheme for coarse levels}  \label{Appendix A}
The stencils of the Laplace and wavenumber operators for interior points on the fourth-, fifth- and sixth-level coarse grid read as

$A_{8h}=\frac{1}{4096} \cdot \frac{1}{1024} \cdot \frac{1}{1024} \cdot \frac{1}{h^2} \cdot$
\[
\resizebox{\textwidth}{!}{$
\begin{bmatrix}\begin{array}{rrrrrrr}
-10395 & -887166 & -7871637 & -15491748 & -7871637 & -887166 & -10395 \\
-887166 & -39105612 & -215169378 & -348459432 & -215169378 & -39105612 & -887166 \\
-7871637 & -215169378 & -265120059 & 413761124 & -265120059 & -215169378 & -7871637 \\
-15491748 & -348459432 & 413761124 & 2809129936 & 413761124 & -348459432 & -15491748 \\
-7871637 & -215169378 & -265120059 & 413761124 & -265120059 & -215169378 & -7871637 \\
-887166 & -39105612 & -215169378 & -348459432 & -215169378 & -39105612 & -887166 \\
-10395 & -887166 & -7871637 & -15491748 & -7871637 & -887166 & -10395 
\end{array} 
\end{bmatrix}
$},
\]

$K_{8h}=\frac{1}{4096} \cdot \frac{1}{4096} \cdot \frac{1}{1024} \cdot k^2 \cdot$
\[
\resizebox{\textwidth}{!}{$
\begin{bmatrix}\begin{array}{rrrrrrr}
27225 & 3939210 & 40768695 & 83544780 & 40768695 & 3939210 & 27225 \\
3939210 & 569967876 & 5898859542 & 12088170168 & 5898859542 & 569967876 & 3939210 \\
40768695 & 5898859542 & 61050008889 & 125106029556 & 61050008889 & 5898859542 & 40768695 \\
83544780 & 12088170168 & 125106029556 & 256372094224 & 125106029556 & 12088170168 & 83544780 \\
40768695 & 5898859542 & 61050008889 & 125106029556 & 61050008889 & 5898859542 & 40768695 \\
3939210 & 569967876 & 5898859542 & 12088170168 & 5898859542 & 569967876 & 3939210 \\
27225 & 3939210 & 40768695 & 83544780 & 40768695 & 3939210 & 27225 
\end{array} 
\end{bmatrix}
$}, 
\]

$A_{16h}=
\frac{1}{4096} \cdot \frac{1}{1024} \cdot \frac{1}{1024} \cdot \frac{1}{1024} \cdot \frac{1}{h^2} \cdot$
\[
\resizebox{\textwidth}{!}{$ 
\begin{bmatrix} \begin{array}{rrrrrrr} 
-13491387 & -1011388446 & -8590720245 & -16705596516 & -8590720245 & -1011388446 & -13491387 \\ -1011388446 & -41427399756 & -220811304386 & -353095695272 & -220811304386 & -41427399756 & -1011388446 \\ -8590720245 & -220811304386 & -262703195227 & 427978620452 & -262703195227 & -220811304386 & -8590720245 \\ -16705596516 & -353095695272 & 427978620452 & 2827174335440 & 427978620452 & -353095695272 & -16705596516 \\ -8590720245 & -220811304386 & -262703195227 & 427978620452 & -262703195227 & -220811304386 & -8590720245 \\ -1011388446 & -41427399756 & -220811304386 & -353095695272 & -220811304386 & -41427399756 & -1011388446 \\ -13491387 & -1011388446 & -8590720245 & -16705596516 & -8590720245 & -1011388446 & -13491387 
\end{array} 
\end{bmatrix} $},
\]

$K_{16h}=\frac{1}{4096} \cdot \frac{1}{4096} \cdot \frac{1}{1024} \cdot \frac{1}{1024} \cdot k^2 \cdot$
\[
\resizebox{\textwidth}{!}{$
\begin{bmatrix}\begin{array}{rrrrrrr}
158684409 & 19907719338 & 199630340247 & 405976871820 & 199630340247 & 19907719338 & 158684409 \\
19907719338 & 2497518765316 & 25044582577654 & 50931743533240 & 25044582577654 & 2497518765316 & 19907719338 \\
199630340247 & 25044582577654 & 251141703197401 & 510732601675060 & 251141703197401 & 25044582577654 & 199630340247 \\
405976871820 & 50931743533240 & 510732601675060 & 1038647851363600 & 510732601675060 & 50931743533240 & 405976871820 \\
199630340247 & 25044582577654 & 251141703197401 & 510732601675060 & 251141703197401 & 25044582577654 & 199630340247 \\
19907719338 & 2497518765316 & 25044582577654 & 50931743533240 & 25044582577654 & 2497518765316 & 19907719338 \\
158684409 & 19907719338 & 199630340247 & 405976871820 & 199630340247 & 19907719338 & 158684409 
\end{array}
\end{bmatrix}
$},
\]

$A_{32h}=\frac{1}{4096} \cdot \frac{1}{1024} \cdot \frac{1}{1024} \cdot \frac{1}{1024} \cdot \frac{1}{1024} \cdot \frac{1}{h^2} \cdot$
\[
\resizebox{\textwidth}{!}{$
\begin{bmatrix}\begin{array}{rrrrrrr}
-14618265915 & -1063059274398 & -8934400311925 & -17322732417892 & -8934400311925 & -1063059274398 & -14618265915 \\
-1063059274398 & -42774103061580 & -226217899925314 & -360608941958056 & -226217899925314 & -42774103061580 & -1063059274398 \\
-8934400311925 & -226217899925314 & -266795319274715 & 439293870677284 & -266795319274715 & -226217899925314 & -8934400311925 \\
-17322732417892 & -360608941958056 & 439293870677284 & 2882610253296592 & 439293870677284 & -360608941958056 & -17322732417892 \\
-8934400311925 & -226217899925314 & -266795319274715 & 439293870677284 & -266795319274715 & -226217899925314 & -8934400311925 \\
-1063059274398 & -42774103061580 & -226217899925314 & -360608941958056 & -226217899925314 & -42774103061580 & -1063059274398 \\
-14618265915 & -1063059274398 & -8934400311925 & -17322732417892 & -8934400311925 & -1063059274398 & -14618265915 
\end{array}
\end{bmatrix}
$},
\]

$K_{32h}=\frac{1}{4096} \cdot \frac{1}{4096} \cdot \frac{1}{1024} \cdot \frac{1}{1024} \cdot \frac{1}{1024} \cdot k^2 \cdot$
\[
\resizebox{\textwidth}{!}{$
\begin{bmatrix}\begin{array}{rrrrrrr}
706549519225 & 85722084590890 & 852977249303575 & 1731387418334860 & 852977249303575 & 85722084590890 & 706549519225 \\
85722084590890 & 10400227566028036 & 103487421517331808 & 210060490730926272 & 103487421517331824 & 10400227566028036 & 85722084590890 \\
852977249303575 & 103487421517331840 & 1029751161146567936 & 2090206046973178368 & 1029751161146568192 & 103487421517331792 & 852977249303575 \\
1731387418334860 & 210060490730926208 & 2090206046973178624 & 4242735025361522688 & 2090206046973178624 & 210060490730926208 & 1731387418334860 \\
852977249303575 & 103487421517331840 & 1029751161146568064 & 2090206046973178624 & 1029751161146567936 & 103487421517331792 & 852977249303575 \\
85722084590890 & 10400227566028036 & 103487421517331808 & 210060490730926208 & 103487421517331808 & 10400227566028036 & 85722084590890 \\
706549519225 & 85722084590890 & 852977249303575 & 1731387418334860 & 852977249303575 & 85722084590890 & 706549519225 
\end{array}
\end{bmatrix}
$}.
\]
 
\section{\hl{Performance analysis using roofline model}}  
\label{Appendix B}
This appendix presents a performance analysis of matrix-vector multiplication operations $\mathbf{v}=\mathbf{A}\mathbf{u}$ comparing our matrix-free implementation with traditional CSR matrix-based approaches. The analysis focuses specifically on the Helmholtz operator with variable wavenumber, using five-point stencil discretization. We consider the matrix-vector multiplication as it constitutes the primary computational kernel in preconditioned Krylov subspace methods, typically accounting for the majority of computational time. Our analysis employs the roofline model \cite{williamsRooflineInsightfulVisual2009}, a performance model that bounds computational kernel performance based on peak computational performance and memory bandwidth limitations. As matrix-vector multiplication is typically memory-bound, we focus on arithmetic intensity ($I$), defined as the ratio of floating-point operations (FLOPs) to memory accesses:  
  
\[
    I = \frac{\text{Total FLOPs}}{\text{Total Bytes Accessed}} \]

\subsection{Analysis of matrix-free implementation}  

The matrix-free implementation directly applies the five-point stencil operation for the discrete Laplacian operator combined with the wavenumber term. For variable wavenumber $\mathbf{k}$, the implementation of matrix-vector multiplication can be expressed in the following computational kernel (in Fortran): 

\begin{lstlisting} 
! Pre-computed stencil coefficients  
ap = -4.d0/h**2  
as = aw = ae = an = 1.d0/h**2  
do j = 1, ny  
    do i = 1, nx   
        v(i,j) = as * u(i,j-1)              & ! South neighbor
                + aw * u(i-1,j)             & ! West neighbor  
                + (ap - k(i,j)**2) * u(i,j) & ! Center point  
                + ae * u(i+1,j)             & ! East neighbor  
                + an * u(i,j+1)               ! North neighbor  
    end do  
end do  
\end{lstlisting} 

In analyzing the memory access patterns, we consider the memory access per grid point operation. The implementation requires reading of five vector elements (center and four neighbors), each consuming 8 bytes in double precision. Additionally, we need to access one wavenumber value (8 bytes) and write one result value (8 bytes). Therefore, the total memory access per point is bounded by 56 bytes.  

The computational intensity involves five multiplications for the stencil coefficients, four additions for combining the stencil components, and two additional operations (one square operation and one subtraction) for the wavenumber term. This results in 11 floating-point operations per grid point.  

Consequently, the arithmetic intensity for the matrix-free implementation is:  
\[I_{MF} \geq \frac{11}{56} \approx 0.1964 \text{ FLOPs/byte} \]

\subsection{Analysis of CSR matrix-based implementation}  

The compressed sparse row (CSR) format represents the sparse matrix $A$ using three arrays: values (\verb|A%value|), column indices (\verb|A%col_indices|), and row pointers (\verb|A%row_ptr|) \cite{saad2003iterative}. The CSR format implementation is structured as follows:

\begin{lstlisting} 
do i = 1, A%nrow  
    v(i) = 0.d0  
    do j = A%row_ptr(i), A%row_ptr(i+1)-1  
        v(i) = v(i) + A%values(j) * u(A%col_indices(j))  
    end do  
end do  
\end{lstlisting}  

The memory access pattern for CSR implementation is more complex. For each non-zero element, we must read the matrix value (8 bytes), the column index (4 bytes), and access the corresponding vector element (8 bytes, assuming the vector is too large to fit into the cache). Additional memory operations include accessing row pointers (4 bytes per row) and reading/writing the result vector (16 bytes per row). For our five-point stencil case, with five nonzero elements per row, the total memory access is bounded by 120 bytes per row.  

The computation for each non-zero element requires one multiplication and one addition, resulting in 10 total FLOPs per row. Thus, the arithmetic intensity for the CSR implementation is:  
\[I_{CSR} \geq \frac{10}{120} \approx 0.0833 \text{ FLOPs/byte}\] 

Based on this theoretical analysis, we expect the matrix-free implementation to outperform the CSR matrix-based implementation by approximately a factor of $2.35$.   

\subsection{Numerical Validation}  
To validate our theoretical analysis, we conducted extensive performance measurements comparing both implementations. For each grid size, we performed 100 consecutive matrix-vector multiplications to obtain statistically stable performance measurements. The performance metrics are reported in billions of floating-point operations per second (GFLOPs/s), calculated using the theoretical operation count for each implementation. 

\begin{table}[htbp]  
\caption{Performance comparison of matrix-free and CSR matrix-based implementations for matrix-vector multiplications.}  
\label{tab:performance}  
\centering  
\scalebox{0.8}{
\begin{tabular}{rccc}  
\hline\noalign{\smallskip}
Problem size   & Matrix-free & CSR-Matrix & Performance \\
$N$    & (GFLOPs/s) & (GFLOPs/s) & Ratio \\
\boldline
289         & 4.5918 & 2.5298 & 1.82 \\
1,089       & 5.1946 & 2.8848 & 1.80 \\
4,225       & 6.1535 & 2.8094 & 2.19 \\
16,641      & 5.0961 & 2.8327 & 1.80 \\
66,049      & 6.1740 & 2.6128 & 2.36 \\
263,169     & 6.1361 & 2.8208 & 2.18 \\
1,050,625   & 6.2861 & 2.3977 & 2.62 \\
4,198,401   & 5.9456 & 1.9992 & 2.97 \\
16,785,409  & 5.5556 & 1.9200 & 2.89 \\
67,125,249  & 5.4626 & 1.8979 & 2.88 \\
268,468,225 & 5.4626 & 1.8958 & 2.88 \\
\hline\noalign{\smallskip} 
\end{tabular}  }
\end{table}  

Table~\ref{tab:performance} presents the experimental results, which strongly support our theoretical analysis. The matrix-free implementation consistently achieves superior performance, with the advantage becoming more pronounced as the problem size increases. For larger problem sizes ($N > 10^6$), we observe performance improvements approaching a factor of $3$, exceeding our theoretical prediction of $2.35$. This enhanced performance can be attributed to memory hierarchy effects. The matrix-free implementation exhibits superior cache utilization, particularly for large-scale problems where the memory access patterns of the CSR format become increasingly inefficient.    

The characteristics of these implementations have significant implications for parallel computing performance. The matrix-free implementation's regular memory access patterns facilitate better parallel efficiency through predictable memory access and reduced NUMA (Non-Uniform Memory Access) effects. Furthermore, in distributed memory systems, the matrix-free approach minimizes communication overhead, requiring only ghost point exchanges along subdomain boundaries. These parallel computing advantages, combined with the superior cache utilization observed in our sequential tests, suggest even more pronounced performance benefits in parallel computing environments, particularly for large-scale problems on distributed memory systems.  

\begin{acknowledgements}
We would like to acknowledge the support of the China Scholarship Council (No. 202006230087).
\end{acknowledgements}

\section*{Data availability}
The source codes supporting the findings of this study are openly available in the GitHub repository \texttt{paraMADP}, DOI \href{https://doi.org/10.5281/zenodo.13143940}{10.5281/zenodo.13143939}. Additional materials or datasets used in this research can be made available upon reasonable request to the corresponding author.

\section*{Conflict of interest}
The authors declare that they have no conflict of interest.

\bibliographystyle{spmpsci}      % mathematics and physical sciences
\bibliography{manuscript}   % name your BibTeX data base

\begin{thebibliography}{10}
\providecommand{\url}[1]{{#1}}
\providecommand{\urlprefix}{URL }
\expandafter\ifx\csname urlstyle\endcsname\relax
  \providecommand{\doi}[1]{DOI~\discretionary{}{}{}#1}\else
  \providecommand{\doi}{DOI~\discretionary{}{}{}\begingroup \urlstyle{rm}\Url}\fi

\bibitem{adrianiAsymptoticSpectralProperties2024}
Adriani, A., Sormani, R.L., Tablino-Possio, C., Krause, R., Serra-Capizzano, S.: Asymptotic spectral properties and preconditioning of an approximated nonlocal {H}elmholtz equation with {C}aputo fractional {L}aplacian and variable coefficient wave number \${$\mu\$$} (2024).
\newblock \urlprefix\url{https://arxiv.org/abs/2402.10569}

\bibitem{babuskaOptimalLocalApproximation2011}
Babuska, I., Lipton, R.: Optimal local approximation spaces for generalized finite element methods with application to multiscale problems.
\newblock Multiscale Model. Simul. \textbf{9}(1), 373--406 (2011).
\newblock \doi{10.1137/100791051}

\bibitem{babuska1997pollution}
Babuska, I.M., Sauter, S.A.: Is the pollution effect of the {FEM} avoidable for the {H}elmholtz equation considering high wave numbers?
\newblock SIAM J. Numer. Anal. \textbf{34}(6), 2392--2423 (1997).
\newblock \doi{10.1137/S0036142994269186}

\bibitem{bootlandComparisonCoarseSpaces2021}
Bootland, N., Dolean, V., Jolivet, P., Tournier, P.H.: A comparison of coarse spaces for {H}elmholtz problems in the high frequency regime.
\newblock Comput. Math. Appl. \textbf{98}, 239--253 (2021).
\newblock \doi{10.1016/j.camwa.2021.07.011}

\bibitem{calandra2013improved}
Calandra, H., Gratton, S., Pinel, X., Vasseur, X.: An improved two-grid preconditioner for the solution of three-dimensional {H}elmholtz problems in heterogeneous media.
\newblock Numer. Linear Algebra Appl. \textbf{20}(4), 663--688 (2013).
\newblock \doi{10.1002/nla.1860}

\bibitem{calandra2017geometric}
Calandra, H., Gratton, S., Vasseur, X.: A geometric multigrid preconditioner for the solution of the {H}elmholtz equation in three-dimensional heterogeneous media on massively parallel computers.
\newblock In: Modern Solvers for {H}elmholtz Problems, pp. 141--155. Springer (2017).
\newblock \doi{10.1007/978-3-319-28832-1_6}

\bibitem{etna_vol59_pp270-294}
Chen, J., Dwarka, V., Vuik, C.: A matrix-free parallel solution method for the three-dimensional heterogeneous {H}elmholtz equation.
\newblock Electron. Trans. Numer. Anal. \textbf{59}, 270--294 (2023).
\newblock \doi{10.1553/etna_vol59s270}

\bibitem{jchen2D2022}
Chen, J., Dwarka, V., Vuik, C.: Matrix-free parallel preconditioned iterative solvers for the 2{D} {H}elmholtz equation discretized with finite differences.
\newblock In: Scientific Computing in Electrical Engineering, pp. 61--68. Springer Nature Switzerland (2024).
\newblock \doi{10.1007/978-3-031-54517-7_7}

\bibitem{chen2023matrixfree2ldef}
Chen, J., Dwarka, V., Vuik, C.: A matrix-free parallel two-level deflation preconditioner for two-dimensional heterogeneous {H}elmholtz problems.
\newblock J. Comput. Phys. p. 113264 (2024).
\newblock \doi{10.1016/j.jcp.2024.113264}

\bibitem{chupengWavenumberExplicitConvergence2023}
Chupeng, M., Alber, C., Scheichl, R.: Wavenumber explicit convergence of a multiscale generalized finite element method for heterogeneous {H}elmholtz problems.
\newblock SIAM J. Numer. Anal. \textbf{61}(3), 1546--1584 (2023).
\newblock \doi{10.1137/21M1466748}

\bibitem{cocquetHowLargeShift2017}
Cocquet, P.H., Gander, M.J.: How large a shift is needed in the shifted {H}elmholtz preconditioner for its effective inversion by multigrid?
\newblock SIAM J. Sci. Comput. \textbf{39}(2), A438--A478 (2017).
\newblock \doi{10.1137/15M102085X}

\bibitem{DHPC2022}
{D}elft {H}igh {P}erformance {C}omputing~{C}entre ({DHPC}): {D}elft{B}lue {S}upercomputer ({P}hase 2).
\newblock \url{https://www.tudelft.nl/dhpc/ark:/44463/DelftBluePhase2} (2024)

\bibitem{drzisga2023semi}
Drzisga, D., K{\"o}ppl, T., Wohlmuth, B.: A semi matrix-free twogrid preconditioner for the {H}elmholtz equation with near optimal shifts.
\newblock J. Sci. Comput. \textbf{95}(3), 82 (2023).
\newblock \doi{10.1007/s10915-023-02195-5}

\bibitem{drzisga2020stencil}
Drzisga, D., R\"{u}de, U., Wohlmuth, B.: Stencil scaling for vector-valued {PDE}s on hybrid grids with applications to generalized {N}ewtonian fluids.
\newblock SIAM J. Sci. Comput. \textbf{42}(6), B1429--B1461 (2020).
\newblock \doi{10.1137/19M1267891}

\bibitem{dwarka2020scalable}
Dwarka, V., Vuik, C.: Scalable convergence using two-level deflation preconditioning for the {H}elmholtz equation.
\newblock SIAM J. Sci. Comput. \textbf{42}(2), A901--A928 (2020).
\newblock \doi{10.1137/18M1192093}

\bibitem{dwarka2020scalablemultilevel}
Dwarka, V., Vuik, C.: Scalable multi-level deflation preconditioning for highly indefinite time-harmonic waves.
\newblock J. Comput. Phys. \textbf{469}, 111327 (2022).
\newblock \doi{10.1016/j.jcp.2022.111327}

\bibitem{Elman2001}
Elman, H.C., Ernst, O.G., O'Leary, D.P.: A multigrid method enhanced by {K}rylov subspace iteration for discrete {H}elmholtz equations.
\newblock SIAM J. Sci. Comput. \textbf{23}(4), 1291--1315 (2001).
\newblock \doi{10.1137/S1064827501357190}

\bibitem{erlangga2005robust}
Erlangga, Y.A.: A robust and efficient iterative method for the numerical solution of the {H}elmholtz equation.
\newblock Ph.D. thesis, Delft University of Technology (2005).
\newblock \urlprefix\url{http://resolver.tudelft.nl/uuid:af9be715-6ebf-4fc1-b948-ebd9d2c4167b}

\bibitem{erlangga2008multilevel}
Erlangga, Y.A., Nabben, R.: Multilevel projection-based nested {K}rylov iteration for boundary value problems.
\newblock SIAM J. Sci. Comput. \textbf{30}(3), 1572--1595 (2008)

\bibitem{erlangga2006novel}
Erlangga, Y.A., Oosterlee, C.W., Vuik, C.: A novel multigrid based preconditioner for heterogeneous {H}elmholtz problems.
\newblock SIAM J. Sci. Comput. \textbf{27}(4), 1471--1492 (2006).
\newblock \doi{10.1137/040615195}

\bibitem{erlangga2004class}
Erlangga, Y.A., Vuik, C., Oosterlee, C.W.: On a class of preconditioners for solving the {H}elmholtz equation.
\newblock Appl. Numer. Math. \textbf{50}(3-4), 409--425 (2004).
\newblock \doi{10.1016/j.apnum.2004.01.009}

\bibitem{gander2019class}
Gander, M.J., Zhang, H.: A class of iterative solvers for the {H}elmholtz equation: Factorizations, sweeping preconditioners, source transfer, single layer potentials, polarized traces, and optimized {S}chwarz methods.
\newblock Siam Review \textbf{61}(1), 3--76 (2019).
\newblock \doi{10.1137/16M109781X}

\bibitem{Gordon2013Robust}
Gordon, D., Gordon, R.: Robust and highly scalable parallel solution of the {H}elmholtz equation with large wave numbers.
\newblock J. Comput. Appl. Math. \textbf{237}(1), 182--196 (2013).
\newblock \doi{10.1016/j.cam.2012.07.024}

\bibitem{grahamDomainDecompositionPreconditioning2017}
Graham, I.G., Spence, E.A., Vainikko, E.: Domain decomposition preconditioning for high-frequency {H}elmholtz problems with absorption.
\newblock Math. Comp. \textbf{86}(307), 2089--2127 (2017).
\newblock \doi{10.1090/mcom/3190}

\bibitem{kim2002multigrid}
Kim, S., Kim, S.: Multigrid simulation for high-frequency solutions of the {H}elmholtz problem in heterogeneous media.
\newblock SIAM J. Sci. Comput. \textbf{24}(2), 684--701 (2002).
\newblock \doi{10.1137/S1064827501385426}

\bibitem{Kononov2007Numerical}
Kononov, A.V., Riyanti, C.D., de~Leeuw, S.W., Oosterlee, C.W., Vuik, C.: Numerical performance of a parallel solution method for a heterogeneous 2{D} {H}elmholtz equation.
\newblock Comput. Vis. Sci. \textbf{11}(3), 139--146 (2007).
\newblock \doi{10.1007/s00791-007-0069-6}

\bibitem{liPreconditioningTechniqueBased2023}
Li, T.Y., Chen, F., Sun, H.W., Sun, T.: Preconditioning technique based on sine transformation for nonlocal {H}elmholtz equations with fractional {L}aplacian.
\newblock J. Sci. Comput. \textbf{97}(1), 17 (2023).
\newblock \doi{10.1007/s10915-023-02332-0}

\bibitem{linAbsolutevalueBasedPreconditioner2024}
Lin, X., Li, C., Hon, S.: Absolute-value based preconditioner for complex-shifted {L}aplacian systems (2024).
\newblock \urlprefix\url{https://arxiv.org/abs/2408.00488}

\bibitem{lu2019robust}
Lu, P., Xu, X.: A robust multilevel preconditioner based on a domain decomposition method for the {H}elmholtz equation.
\newblock J. Sci. Comput. \textbf{81}, 291--311 (2019).
\newblock \doi{10.1007/s10915-019-01015-z}

\bibitem{maNovelDesignAnalysis2022}
Ma, C., Scheichl, R., Dodwell, T.: Novel design and analysis of generalized finite element methods based on locally optimal spectral approximations.
\newblock SIAM J. Numer. Anal. \textbf{60}(1), 244--273 (2022).
\newblock \doi{10.1137/21M1406179}

\bibitem{plessix2003separation}
Plessix, R.E., Mulder, W.A.: Separation-of-variables as a preconditioner for an iterative {H}elmholtz solver.
\newblock Appl. Numer. Math. \textbf{44}(3), 385--400 (2003).
\newblock \doi{10.1016/S0168-9274(02)00165-4}

\bibitem{pulchHelmholtzEquationUncertainties2024}
Pulch, R., S{\`e}te, O.: The {H}elmholtz equation with uncertainties in the wavenumber.
\newblock J. Sci. Comput. \textbf{98}(3), 60 (2024).
\newblock \doi{10.1007/s10915-024-02450-3}

\bibitem{riyanti2007parallel}
Riyanti, C., Kononov, A., Erlangga, Y., Vuik, C., Oosterlee, C., Plessix, R.E., Mulder, W.: A parallel multigrid-based preconditioner for the 3{D} heterogeneous high-frequency {H}elmholtz equation.
\newblock J. Comput. Phys. \textbf{224}(1), 431--448 (2007).
\newblock \doi{10.1016/j.jcp.2007.03.033}

\bibitem{saad2003iterative}
Saad, Y.: Iterative Methods for Sparse Linear Systems, second edn.
\newblock Society for Industrial and Applied Mathematics (2003).
\newblock \doi{10.1137/1.9780898718003}

\bibitem{sheikh2016accelerating}
Sheikh, A.H., Lahaye, D., Ramos, L.G., Nabben, R., Vuik, C.: Accelerating the shifted {L}aplace preconditioner for the {H}elmholtz equation by multilevel deflation.
\newblock J. Comput. Phys. \textbf{322}, 473--490 (2016).
\newblock \doi{10.1016/j.jcp.2016.06.025}

\bibitem{sheikh2013convergence}
Sheikh, A.H., Lahaye, D., Vuik, C.: On the convergence of shifted laplace preconditioner combined with multilevel deflation.
\newblock Numer. Linear Algebra Appl. \textbf{20}(4), 645--662 (2013).
\newblock \doi{10.1002/nla.1882}

\bibitem{sourbierThreedimensionalParallelFrequencydomain2011}
Sourbier, F., Haidar, A., Giraud, L., Ben-Hadj-Ali, H., Operto, S., Virieux, J.: Three-dimensional parallel frequency-domain visco-acoustic wave modelling based on a hybrid direct/iterative solver.
\newblock Geophysical Prospecting \textbf{59}(5), 834--856 (2011).
\newblock \doi{10.1111/j.1365-2478.2011.00966.x}

\bibitem{tang2009comparison}
Tang, J.M., Nabben, R., Vuik, C., Erlangga, Y.A.: Comparison of two-level preconditioners derived from deflation, domain decomposition and multigrid methods.
\newblock J. Sci. Comput. \textbf{39}, 340--370 (2009).
\newblock \doi{10.1007/s10915-009-9272-6}

\bibitem{taus2020sweeps}
Taus, M., Zepeda-N{\'u}{\~n}ez, L., Hewett, R.J., Demanet, L.: L-sweeps: A scalable, parallel preconditioner for the high-frequency {H}elmholtz equation.
\newblock J. Comput. Phys. \textbf{420}, 109706 (2020).
\newblock \doi{10.1016/j.jcp.2020.109706}

\bibitem{tournierNumericalModelingHighspeed2017}
Tournier, P.H., Bonazzoli, M., Dolean, V., Rapetti, F., Hecht, F., Nataf, F., Aliferis, I., El~Kanfoud, I., Migliaccio, C., De~Buhan, M., Darbas, M., Semenov, S., Pichot, C.: Numerical modeling and high-speed parallel computing: {{New}} perspectives on tomographic microwave imaging for brain stroke detection and monitoring.
\newblock IEEE Antennas Propag. Mag. \textbf{59}(5), 98--110 (2017).
\newblock \doi{10.1109/MAP.2017.2731199}

\bibitem{Versteeg_1991_ME}
Versteeg, R.: The marmousi experience: Velocity model determination on a synthetic complex data set.
\newblock The Leading Edge \textbf{13}(9), 927--936 (1994)

\bibitem{williamsRooflineInsightfulVisual2009}
Williams, S., Waterman, A., Patterson, D.: Roofline: {{An}} insightful visual performance model for multicore architectures.
\newblock Commun. ACM \textbf{52}(4), 65--76 (2009).
\newblock \doi{10.1145/1498765.1498785}

\bibitem{yovelLFAtunedMatrixfreeMultigrid2024}
Yovel, R., Treister, E.: {{LFA-tuned}} matrix-free multigrid method for the elastic {H}elmholtz equation.
\newblock SIAM J. Sci. Comput. pp. S1--S21 (2024).
\newblock \doi{10.1137/23M1583466}

\end{thebibliography}


\end{document}